\documentclass[11pt]{article}
\usepackage{amsmath,amsfonts,amssymb}
\usepackage{multirow} 
\usepackage{float}

\usepackage[section]{placeins} 

\usepackage{graphicx}
\usepackage{epstopdf}
\usepackage{bm}
\usepackage{color}
\usepackage{lscape}
\usepackage{url}
\usepackage{hyperref}

\DeclareMathOperator\erfc{erfc}
\DeclareMathOperator\erf{erf}

\newcommand{\bd}[1]{\mathbf{#1}}
\newcommand{\smint}{\textstyle{\int}}
\newcommand{\tb}{\mathbf}
\newcommand{\ff}{{\tb f}}
\newcommand{\qq}{{\tb q}}
\newcommand{\uu}{{\tb u}}
\newcommand{\xx}{{\tb x}}
\newcommand{\yy}{{\tb y}}
\newcommand{\zz}{{\tb z}}
\newcommand{\rr}{{\tb r}}
\newcommand{\nn}{{\tb n}}
\newcommand{\kk}{{\tb k}}
\newcommand{\fftw}{{\tilde \ff}}
\newcommand{\xh}{{\hat x}}
\newcommand{\xxh}{{\hat \xx}}
\newcommand{\al}{\alpha}
\newcommand{\del}{\delta}
\newcommand{\eps}{\epsilon}
\newcommand{\lam}{\lambda}
\newcommand{\pa}{\partial}
\newcommand{\beq}{\begin{equation}}
\newcommand{\eeq}{\end{equation}}
\newcommand{\Sl}{\mathcal{S}}
\newcommand{\Dl}{\mathcal{D}}

\definecolor{red}{rgb}{1.0,0.0,0.0}

\title{
Extrapolated Regularization of \\ Nearly Singular Integrals on Surfaces}

\author{J. Thomas Beale\thanks{Department of Mathematics, Duke University, Durham, NC, 27708 USA beale@math.duke.edu} 
\and Svetlana Tlupova\thanks{Department of Mathematics,
Farmingdale State College, SUNY, Farmingdale, NY 11735, USA tlupovs@farmingdale.edu}}

\date{\today}

\begin{document}

\maketitle

\begin{abstract}
We present a method for computing nearly singular integrals that occur when single or double layer surface integrals, for harmonic potentials or Stokes flow, are evaluated at points nearby.  Such values could be needed in solving an integral equation when one surface is close to another or to obtain values at grid points.  We replace the singular kernel with a regularized version having a length parameter $\del$ in order to control discretization error.  
Analysis near the singularity leads to an expression for the error due to regularization which has terms with unknown coefficients multiplying known quantities.  By computing the integral with three choices of $\del$
we can solve for an extrapolated value that has regularization error reduced to $O(\del^5)$,
uniformly for target points on or near the surface.
  In examples with
$\del/h$ constant and moderate resolution we observe total error about $O(h^5)$ close to the surface.  For convergence as $h \to 0$ we can choose $\del$ proportional to $h^q$ with $q < 1$ to ensure the discretization error is dominated by the regularization error.  With $q = 4/5$ we find errors about $O(h^4)$.  For harmonic potentials we extend the approach to a version with $O(\del^7)$ regularization; it typically has smaller errors but the order of accuracy is less predictable.
\end{abstract}



\section{Introduction}
The evaluation of singular or nearly singular surface integrals, on or near the surface, requires special care.  Here we are concerned with single and double layer integrals for harmonic potentials or for Stokes flow.  One of several possible approaches is to regularize the singular kernel in order to control the discretization error.  A natural choice is to replace the $1/r$ singularity in the single layer potential with $\erf(r/\delta)/r$, where $\erf$ is the error function and $\delta$ is a numerical parameter setting the length scale of the regularization.  This replacement introduces an additional error due to smoothing.  For the singular case, evaluating at points on the surface, we can modify the choice of regularization so that the new error is $O(\delta^5)$; see \cite{b04,byw,tbjcp}.  The nearly singular case, evaluation at points near the surface, could be needed e.g. 
in solving integral equations
when surfaces are close together or to obtain values at grid points.
For this case, in the previous work, we used
analysis near the singularity to derive corrections which leave a remaining error of $O(\delta^3)$.  It does not seem practical to extend the corrections to higher order. In the present work we show by local analysis that the simpler regularization can be used with extrapolation, rather than corrections, to improve the error to $O(\delta^5)$ in the nearly singular case.  For $\yy$ on or near the surface, at signed distance $b$,
if $\Sl$ is the single layer potential with some density function and $\Sl_\del$ is the regularized integral, we show that
\beq \label{front}  \Sl_\del(\yy) =
\Sl(\yy) +  C_1\del I_0(b/\del)  + C_2\del^3 I_2(b/\del)
          + O(\del^5)  \eeq
uniformly for $\yy$ as $\delta \to 0$,
where $I_0$ and $I_2$ are certain integrals, known explicitly, and $C_1$, $C_2$ are coefficients which depend on $\yy$, $b$, the surface, and the density function.  We can regard $\Sl$, $C_1$, $C_2$ as unknowns at one point $\yy$.  Our strategy is to calculate the
regularized integrals $\Sl_\del$
 for three different choices of $\del$ and then solve for $\Sl$, within $O(\del^5$), from the system of three equations. 
We treat the double layer potential in a similar way, as well as the single and double layer  integrals for Stokes flow.
We comment on the Helmholtz equation.
For the harmonic potentials we extend the approach to a method with $O(\del^7)$ regularization error; it requires four choices of $\del$ rather than three.

To compute the integrals we use a quadrature rule for surface integrals for which the quadrature points are points where the surface intersects lines in a three-dimensional grid and the weights are determined by the normal vector to the surface.  
It is high order accurate for smooth integrands; for the nearly singular integrals the accuracy depends on $\delta$ as well as the grid size $h$.
The regularization enables us to make the integrand smooth enough to discretize without special treatment near the singularity.
Other quadrature methods could be used if desired.
The total error consists of the regularization error and the error due to discretization.  The  discretization error is low order as $h \to 0$ if $\delta/h$ is fixed,
but it rapidly improves as $\delta/h$ increases; this is explained in Sect.~4.  In our experiments with
$\delta/h$ constant, we typically observe errors about $O(h^5)$ near the surface with moderate resolution, i.e. $h$ not too small, indicating that the regularization error is dominant.  However this trend cannot continue as $h \to 0$.
For rapid convergence as $h \to 0$ we need to increase $\del/h$ to ensure that the discretization error is dominated by the regularization error.  To do this we choose $\del$ proportional to $h^q$, e.g. with
$q = 4/5$, resulting in an error about $O(h^4)$.  
To test the uniform convergence we measure errors at grid points 
within distance $h$ from the surface.
With the fifth order regularization
we see the predicted orders, while for the seventh order method
we typically see smaller errors but the order in $h$ is less predictable,
presumably because of discretization error.

Considerable work has been devoted to the computation of singular integrals such as layer potentials.  Only a portion of this work has concerned nearly singular integrals on surfaces.  Often
values close to the surface are obtained by extrapolating from values further away \cite{ying-biros-zorin-06}, sometimes as part of the quadrature by expansion
(QBX) or hedgehog methods
 \cite{baggetorn,klinttorn,
klockner-barnett-greengard-oneil-13,zorincplx,siegel-tornberg-18}.  In \cite{steinbarn} sources are placed on the opposite side of the surface to produce a kernel independent method.
With the singularity subtraction technique \cite{helsing-13} a most singular part is evaluated analytically leaving a more regular remainder.  In \cite{nitsche}, for the nearly singular axisymmetric case,
the error in computing the most singular part provides a correction. In \cite{perez} an approximation to the density function is used to reduce the singularity.
Regularization has been used extensively to model Stokes flow in biology \cite{cortez2d,cortez3d}; see also \cite{tbjcp}.  Richardson extrapolation has been used for
Stokes flow \cite{gallsmith}.

With Ewald splitting \cite{ew1},\cite{ew2},\cite{ew3},\cite{ew4},\cite{ew5}
 the kernel is written as a localized singular part plus a smooth part so that the two parts can be computed by different methods.  Regularization as used
here could be thought of as a limit case which reduces the singular part so that it becomes a correction, as in \cite{b04,byw,tbjcp}
or treated as an error in the present case.  Integrals for
the heat equation were treated in this way in \cite{greenstrain}, with the history treated as a smooth part.  There is an analogy between the present method and QBX.
In the latter, the value
at a specified point near the boundary is extrapolated from values at points further away along a normal line; increasing the distance is a kind of smoothing, analogous
to the regularization here.  However the two techniques for making the integral smoother are different in practice.

While the choice of numerical method depends on context, the present approach is simple and direct.  The work required is similar to that for a surface integral with smooth integrand, except that three (or four) related integrals must be computed rather than one.  No special gridding or separate treatment of the singularity is needed.  The surface must be moderately smooth, without corners or edges.  Geometric information about the surface is not needed other than normal vectors; further geometry was needed for the corrections of \cite{b04,byw,tbjcp} and in some other methods.  It would be enough for the surface to be known through values of a level set function at grid points nearby.  For
efficiency fast summation methods suitable for regularized kernels \cite{yingradial,rbfstokes,wang-krasny-tlupova-20} could be used.  The approach here is general enough that it should apply to other singular kernels; however, a limitation is discussed at the end of the next section.

Results are described more specifically in Sect.~2.  The analysis leading to \eqref{front} is carried out in Sect.~3.  In Sect.~4 we discuss the quadrature rule and the discretization error.  In Sect.~5 we present numerical examples which illustrate the behavior of the method.  In Sect.~6 we prove that the system of three equations of the form \eqref{front} is solvable, and Sect.~7 has a brief conclusion.

\section
{Summary of results}

For a single layer potential 
\beq \label{sgllayer}
  \Sl(\yy) = \int_\Gamma G(\xx - \yy)f(\xx)\,dS(\xx)\,, \quad
  G(\rr) =  -\frac{1}{4\pi|\rr|} \eeq
on a closed surface $\Gamma$, with given density function $f$, we define the regularized version
\beq \label{sglreg}
   \Sl_\del(\yy) = \int_\Gamma G_\del(\xx - \yy)f(\xx)\,dS(\xx)\,, \quad     G_\delta(\rr) = G(\rr)s_1(|\rr|/\del) \eeq
with 
\beq \label{sglshape}
s_1(r) = \erf(r) = \frac{2}{\sqrt{\pi}}\int_0^r e^{-s^2}\,ds \eeq
Then $G_\delta$ is smooth, with $G_\del(0) = -1/(2\pi^{3/2}\del)$, and
$\erf(r/\del) \to 1$ rapidly as $r/\del$ increases.
Typically $\Sl_\del - \Sl = O(\del)$.
If $\yy$ is near the surface, then $\yy = \xx_0 + b\nn$, where $\xx_0$ is the closest point on $\Gamma$, $\nn$ is the outward normal vector at $\xx_0$, and $b$ is the signed distance.  From a series expansion for $\xx$ near $\xx_0$ and $b$ near $0$ we show in Sect. 3 that
\beq \label{basic}
  \Sl(\yy) +  C_1\del I_0(\lam)  + C_2\del^3 I_2(\lam)
         = \Sl_\del(\yy)  + O(\del^5)  \eeq 
uniformly for $\yy$ near the surface,
where $\lam = b/\del$; $C_1$, $C_2$ are unknown coefficients; and $I_0$ and $I_2$ are integrals occurring in the derivation that are found to be
\beq \label{eye0}
I_0(\lambda) = e^{-\lambda^2}/\sqrt{\pi}
                         - |\lambda|\erfc|\lambda|  \eeq
\beq \label{eye2}
I_2(\lambda) = \frac23\left( 
    (\frac12 - \lambda^2) e^{-\lambda^2}/\sqrt{\pi}
       + |\lambda|^3 \erfc|\lambda| \right) \eeq
Here $\erfc = 1 - \erf$.
To obtain an accurate value of $\Sl$, we calculate the
regularized integrals $\Sl_\del$ for three different choices of $\del$,
at the same $\yy$ with the same grid size $h$,
resulting in a system of three equations with three unknowns.  We can then solve for the exact integral $\Sl$ within error $O(\del^5)$.
We typically choose $\delta_i = \rho_i h$ with $\rho_i = 2, 3, 4$ or $3, 4, 5$. 

To improve the conditioning we write three versions of \eqref{basic} in terms of $\rho$ rather than $\del$,
\beq \label{form5}
   \Sl(\yy) +  c_1\rho_i I_0(\lam_i)  + c_2\rho_i^3 I_2(\lam_i)
         = \Sl_{\del_i}(\yy)  + O(\del_i^5) \,,\quad i = 1,2,3  \eeq
with $\lam_i = b/\del_i$.  It is important that $c_1, c_2$ do not depend on $\del$ or $\lam$.  We solve this $3\times 3$ system for $S$.
The $i$th row is $[1\;,\rho_iI_0(\lam_i)\;\,\rho_i^3I_2(\lam_i)]$; the entries depend only on
$\lam_i$ as well as $\rho_i$.  The value obtained for $\Sl$ has the form
\beq \label{comb}
  \Sl(\yy) = \sum_{i=1}^3 a_iS_{\del_i}(\yy) \,, \quad a_i = a_i(\lam_1,\lam_2,\lam_3)    \eeq
In each case $a_1 + a_2  + a_3 = 1$.
For $\yy \in \Gamma$, $\lam_i = 0$ and
$a_1 = 14/3$, $a_2 = -16/3$, $a_3 = 5/3$ provided $\rho_i = 2, 3, 4$.  As $b$ increases, the coefficients approach
$(1,0,0)$, allowing a gradual transition to the region far enough from $\Gamma$ to omit the regularization.  
It is not obvious that the system \eqref{form5} is solvable, i.e. that the matrix is invertible.  In Sect. 6 we
prove the solvability for any distinct choices of the $\rho_i$.
To ensure the smoothing error is dominant as $h \to 0$ we may choose $\delta = \rho h^q$ with $q < 1$, rather than $q = 1$, to obtain convergence $O(h^{5q})$; see Sect. 4.

For the double layer potential
\beq \label{dbllayer}
  \Dl(\yy) = \int_\Gamma \frac{\pa G(\xx-\yy)}{\pa\nn(\xx)}
        g(\xx)\,dS(\xx) \eeq
the treatment is similar after a subtraction.  Using Green's identities we rewrite \eqref{dbllayer} as
\beq \label{dblsub} 
  \Dl(\yy) = \int_\Gamma \frac{\pa G(\xx-\yy)}{\pa\nn(\xx)}
 [g(\xx) - g(\xx_0)]\,dS(\xx) + \chi(\yy)g(\xx_0) \eeq
where again $\xx_0$ is the closest point on $\Gamma$ and 
$\chi = 1$ for $\yy$ inside, $\chi = 0$ for $\yy$ outside, and
$\chi = \frac12$ on $\Gamma$.  To regularize we replace $\nabla G$
with the gradient of the smooth function $G_\del$, obtaining
\beq \label{gradreg}
  \nabla G_\del(\rr) = \nabla G(\rr)s_2(|\rr|/\del) = 
    \frac{\rr}{4\pi|\rr|^3}s_2(|\rr|/\del)  \eeq
with
\beq \label{dblshape}
     s_2(r) = \erf(r) - (2/\sqrt{\pi})re^{-r^2}  \eeq
Thus
\beq \label{dblreg}
  \Dl_\del(\yy) = \int_\Gamma 
  \frac{\rr\cdot\nn(\xx)}{4\pi|\rr|^3}s_2(|\rr|/\del)
  [g(\xx) - g(\xx_0)]\,dS(\xx) + \chi(\yy)g(\xx_0)\,,
             \quad    \rr = \xx - \yy  \eeq
The expansion for $\Dl_\del - \Dl$ near $\xx_0$ is somewhat different but coincidentally leads to the same relation as in \eqref{form5} with
$\Sl$ and $\Sl_\del$ replaced by $\Dl$ and $\Dl_\del$.  Thus we can solve
for $\Dl$ to $O(\del^5)$ in the same way as for $\Sl$.

There is a straightforward extension to a method with $O(\del^7)$ regularization error.  In equation \eqref{basic} there is now
an additional term $C_3\del^5 I_4(\lambda)$.  There are four unknowns, so that four choices of $\del$ are needed.
Otherwise this version is similar to the original one.  On the other hand, we could use only two choices of
$\del$, omitting the $\del^3$ term in \eqref{basic}, obtaining a version with error $O(\del^3)$.

The special case of evaluation at points $\yy$ on the surface
is important because it is used to solve integral equations for
problems such as the Dirichlet or Neumann problem for harmonic functions.
We could use the procedure described with $b = 0$ and $\lam = 0$.  However in this case we can modify the regularization to obtain $O(\delta^5)$ error more directly \cite{b04,byw}.
For the single layer integral, in place of \eqref{sglreg} we use
\beq \label{sglon}
  G_\del(\rr) = -\frac{s_1^\sharp(|\rr|/\del)}{4\pi|\rr|}\,, \quad
  s_1^\sharp(r) = \erf(r) + \frac{2}{3\sqrt{\pi}}(5r - 2r^3)e^{-r^2} \eeq
For the double layer
we use \eqref{dblreg} with $\chi = \frac12$ and \eqref{dblshape}
replaced by 
\beq \label{dblonshape}
  s_2^\sharp(r) = \erf(r) - \frac{2}{\sqrt{\pi}}
  \left(r - \frac{2r^3}{3}\right)e^{-r^2}  \eeq
We typically use $\del = 3h$ with these formulas for evaluation on the surface \cite{byw,tbjcp}.  They
were derived by imposing conditions to eliminate the leading error \cite{b04}, and the error can be checked using the analysis in the next section.  Formulas with $O(\delta^7)$ error could be produced with the same approach.

The equations of Stokes flow represent the motion of incompressible fluid in the limit of zero Reynolds number; e.g. see \cite{pozbook}.  In the simplest form they are
\beq \label{stokeseqns}
  \Delta\uu - \nabla p = 0\,,\quad \nabla\cdot\uu = 0 \eeq
where $\uu$ is the fluid velocity and $p$ is the pressure.  The primary fundamental solutions for the velocity are the Stokeslet and stresslet,
\begin{subequations}
\begin{align}
	\label{Stokeslet}
	S_{ij}(\bd{y,x}) &= \frac{\del_{ij}}{|\bd{y} - \bd{x}|} + \frac{(y_i - x_i
)(y_j - x_j)}{|\bd{y} - \bd{x}|^3} \\[6pt]
	\label{stresslet}
	T_{ijk} (\bd{y,x}) &= -\frac{6(y_i - x_i)(y_j - x_j)(y_k - x_k)}{|\bd{y} -
 \bd{x}|^5}
\end{align}
\end{subequations}
where $\del_{ij}$ is the Kronecker delta and $i,j,k = 1,2,3$.  They are the kernels for the single and double layer integrals
\begin{subequations}
\begin{align}
	\label{SingleLayer}
	u_i(\bd{y}) &= \frac{1}{8\pi}\int_\Gamma S_{ij}(\bd{y,x}) f_j(\bd{x})
dS(\bd{x}) \\[6pt]
	\label{DoubleLayer}
	v_i(\bd{y}) &= \frac{1}{8\pi}\int_\Gamma T_{ijk} (\bd{y,x}) q_j(\bd{x}) n_k(\bd{x})dS(\bd{x})
\end{align}
\end{subequations}
where $f_j$ and $q_j$ are components of vector quantities $\ff$ and $\qq$ on the surface and $n_k$ is a component of the normal vector $\nn$.  A subtraction can be used in both cases; e.g., see \cite{pozbook},
Sect. 6.4.   With $\xx_0$ as before we rewrite
\eqref{SingleLayer} as
\beq \label{stosglsub}
u_i(\bd{y}) = \frac{1}{8\pi}\int_\Gamma S_{ij}(\bd{y,x})
 [f_j(\xx) - f_k(\xx_0)n_k(\xx_0)n_j(\xx)]\,dS(\bd{x}) \eeq
The subtracted form of \eqref{DoubleLayer} is
\beq \label{stodblsub}
	v_i(\bd{y}) = \frac{1}{8\pi}\int_\Gamma T_{ijk} (\bd{y,x}) [q_j(\bd{x
}) - q_j(\bd{x}_0)] n_k(\bd{x}) dS(\bd{x}) + \chi (\bd{y}) q_i(\bd{x
}_0)
\eeq
To compute \eqref{stosglsub} we replace $S_{ij}$ with the regularized version
\beq \label{Sreg} S_{ij}^\del(\bd{y,x}) = \frac{\del_{ij}}{r}s_1(r/\del)
     + \frac{(y_i - x_i)(y_j - x_j)}{r^3}s_2(r/\del)\,,
         \quad r = |\yy - \xx| \eeq
with $s_1$ and $s_2$ as in \eqref{sglshape},\eqref{dblshape}, resulting in a smooth kernel.

For the Stokes double layer integral we need to rewrite the kernel
so that it will be compatible with the analysis of Sect. 3; see the last paragraph of this section
for further discussion.
For $\yy$ near the surface we have
$\yy = \xx_0 + b\nn_0$ with $\xx_0 \in \Gamma$ and $\nn_0 = \nn(\xx_0)$.
In $T_{ijk}$ we substitute $y_i - x_i = bn_i - \xh_i$ where $n_i$ and $\xh_i$ are the $i$th components of $\nn_0$ and $\xxh = \xx - \xx_0$
and similarly for $j$ and $k$.
The product $(y_i - x_i)(y_j - x_j)(y_k - x_k)$ becomes a sum.  We need to avoid terms in the kernel such as
$b^3/r^5$ or $b^2x_i/r^5$, with $r = |\yy - \xx|$.  To do this
we replace $b^2/r^2$ with  $1 - (r^2 - b^2)/r^2$ 
to introduce factors in the numerator which vanish
at $\xx_0$.  We obtain 
\beq \label{Tsplit}
  T_{ijk} =  T_{ijk}^{(1)} + T_{ijk}^{(2)} =
   - 6\left(\frac{t_{ijk}^{(1)}}{r^3} +
\frac{t_{ijk}^{(2)}
  - (r^2 - b^2)t_{ijk}^{(1)}}{r^5}\right) \eeq
where 
\beq t_{ijk}^{(1)}  = bn_in_jn_k - 
        (\xh_in_jn_k + n_i\xh_jn_k + n_in_j\xh_k)\eeq
\beq t_{ijk}^{(2)} = 
 b(\xh_i\xh_jn_k + \xh_in_j\xh_k + n_i\xh_j\xh_k) - \xh_i\xh_j\xh_k \eeq
and we substitute $r^2 - b^2 = |\xxh|^2 - 2b\xxh\cdot\nn_0$. 
We compute \eqref{stodblsub} with $T_{ijk}$ replaced with the regularized version of \eqref{Tsplit}
\beq \label{Treg} 
T_{ijk}^\del = T_{ijk}^{(1)}s_2(r/\del) + T_{ijk}^{(2)}s_3(r/\del) \eeq
where
\beq \label{s3}
     s_3(r) = \erf(r) 
            - \frac{2}{\sqrt{\pi}}\left(\frac23 r^3 + r\right)e^{-r^2}  \eeq

For both Stokes integrals, calculated in the manner described, we find in Sect. 3 that the error has a form equivalent to \eqref{form5}, and we extrapolate with three choices of
$\del$ as before.  Again for the special case of evaluation on the surface we can obtain an $O(\del^5)$ regularization directly.  Formulas were given in \cite{tbjcp} and an improved formula for the stresslet case was given in \cite{novel}.

A strategy similar to that for the Laplacian
could be used for single or double layer integrals for the Helmholtz equation,
$\Delta u + k^2 u = 0$, which describes waves of a definite frequency.
The usual fundamental solution is $G = -e^{ikr}/4\pi r$.  We could regularize the most singular part,
$-1/4\pi r$ or $-(1 - k^2r^2/2)/4\pi r$ or $-\cos{kr}/4\pi r$, multiplying by $\erf(r/\delta)$,
 and extrapolate as for the Laplacian.  We would not modify the remaining part of $G$.
For the double layer potential we need to use a subtraction again.  We could
do this using a plane wave and Green's third identity (e.g. see \cite{nedelec} Thm. 3.1.1) as has been done before
(e.g. see \cite{planewave}).  We choose a vector
 $\kk$ so that $|\kk| = k$ and, for convenience, $\kk\cdot\nn(\xx_0) = 0$.
 With $\xx_0$ and $\chi(\yy)$ as in \eqref{dblsub} we rewrite the double layer potential as
\begin{multline}\label{Hdbl} 
      \int \frac{\pa G(\xx-\yy)}{\pa\nn(\xx)}g(\xx)\,dS(\xx) = 
  \int \frac{\pa G(\xx-\yy)}{\pa\nn(\xx)}\left[g(\xx) - e^{i\kk\cdot(\xx-\xx_0)}g(\xx_0)\right
]\,dS(\xx) \\
        + g(\xx_0)\int i\kk\cdot\nn(\xx)\,G(\xx-\yy)\,e^{i\kk\cdot(\xx-\xx_0)}\,dS(\xx)
            + \chi(\yy) e^{i\kk\cdot(\yy-\xx_0)}g(\xx_0)   
\end{multline}
 If we regularize only the $1/r$ term we could instead use
\eqref{dblsub} for that part alone.

It appears this method would not be successful if applied directly to the double layer potential or the Stokeslet integral without the subtraction.  There would be 
a term in the integrand proportional to 
$1/r^3$.  The equation \eqref{basic} for the regularization error would then have an additional term which, to first approximation, does not change as $\del$ is varied.  As a result the extrapolated value of the integral becomes unstable as $b \to 0$; i.e.,
the coefficients in the linear combination replacing \eqref{comb} become large as $b \to 0$.  A similar consideration motivates the expression for $T_{ijk}$ above.  For other kernels general techniques to reduce the singularity could be used if necessary, e.g. \cite{perez}.

\section{Local analysis near the singularity}

We derive an expansion for the error due to regularizing a singular integral, when evaluated at a point $\yy$ near the surface $\Gamma$.
The error is uniform with respect to $\yy$.
  The expression obtained leads to the formula \eqref{basic} and the extrapolation strategy used here.  The first few terms of the expansion were used in \cite{b04,byw,tbjcp} to find corrections to $O(\del^3)$.

We begin with the single layer potential \eqref{sgllayer}. The error $\eps$ is the difference between \eqref{sglreg} and \eqref{sgllayer}.
 Given $\yy$ near $\Gamma$,
we assume for convenience that
the closest point on $\Gamma$ is $\xx = 0$.  Then $\yy = b\nn_0$, where
$\nn_0$ is the outward normal at $\xx = 0$ and $b$ is the signed distance from the surface.  We choose coordinates $\al = (\al_1,\al_2)$ on $\Gamma$ near $\xx = 0$ so that $\xx(0) = 0$, the metric tensor $g_{ij} = \del_{ij}$ at $\al = 0$, and the second derivatives $\xx_{ij}$ are normal at $\al = 0$.  E.g., if the tangent plane at $\xx = 0$ is $\{x_3 = 0\}$, we could use $(\al_1,\al_2) = (x_1,x_2)$.  Since the error in the integral is negligible for $\xx$ away from $0$ we can assume the density $f$ is zero outside this coordinate patch, regard it as a function of $\al$, and write the regularization error as
\beq \eps = 
  \int [G_\del(\xx(\al) - \yy) - G(\xx - \yy)]f(\al)\,dS(\al) \eeq
Then
\beq \eps = \frac{1}{4\pi}\int\frac{\erfc(r/\del)}{r}f(\al)\,dS(\al)
\,, \qquad r = |\xx(\al) - \yy|  \eeq

We can expand $\xx$ near $0$ as
\beq \xx(\al) = {\tb T}_1(0)\al_1 + {\tb T}_2(0)\al_2 
 + \sum_{2\leq |\nu|\leq 4} c_\nu \al^\nu D^\nu\xx(0) + O(|\al|^5)  \eeq
Here ${\tb T}_i = \pa\xx/\pa\al_i$, the tangent vector at 
$\xx(\al)$, and we use multi-index notation:  $\nu = (\nu_1,\nu_2)$, $\al^\nu = \al_1^{\nu_1}\al_2^{\nu_2}$, $D^\nu$ is mixed partial derivative of order $(\nu_1,\nu_2)$, and $|\nu| =
\nu_1 + \nu_2$.  We will use the notation $c_\nu$ for generic constants whose value will not be needed.
We first get an expression for  $r^2$.  We start with
\beq |\xx(\al)|^2 = \al_1^2 + \al_2^2 + 
     \sum_{|\nu|=4,5} c_\nu\al^\nu + O(|\al|^6) \eeq
There is no term with $|\nu| = 3$ since the first and second order terms in $\xx$ are orthogonal.  Also
\beq \xx(\al)\cdot \nn_0 = \sum_{2\leq\nu\leq 4}c'_\nu\al^\nu
      + O(|\al|^5) \eeq
Then 
\beq \label{rsq} r^2 = |x(\al) - b\nn_0|^2 = 
  |\al|^2 + b^2 +  \sum_{|\nu|=4,5} c_\nu\al^\nu 
  + b\sum_{2\leq|\nu|\leq 4} c'_\nu\al^\nu
  + O(|\al|^6 + |b||\al|^5) \eeq
  We assume $\Gamma$ is smooth, so that the error terms are uniform with respect to the
  location.

We will make a change of variables $\al = (\al_1,\al_2) \to \xi = (\xi_1,\xi_2)$ defined by
\beq \label{defxi}
|\xi|^2 + b^2 = r^2\,,\quad \xi/|\xi| = \al/|\al| \eeq
This allows us to write the error as
\beq \label{epsxib}
\eps = \frac{1}{4\pi}\int
  \frac{\erfc(\sqrt{|\xi|^2 + b^2}/\del)}{\sqrt{|\xi|^2 + b^2}}
   w(\xi,b)\,d\xi  \eeq
where
\beq \label{wsgl}
 w(\xi,b) = f(\al)\left|\frac{\pa\al}{\pa\xi}\right| |T_1\times T_2|   \eeq
 The expression \eqref{epsxib} will enable us to expand the error in the form we need.  An estimate of
 \eqref{epsxib}, bounding $w$ by a constant, shows that $\epsilon$ decays faster than $e^{-b^2/\del^2}$ and so is negligible for $|b|$ larger
 than $O(\del)$.  Thus we can regard $|b|$ as being at most $O(\delta)$.
  
The mapping $\xi = \xi(\al)$ is close to the identity but it is not smooth at $\al = 0$, so that we cannot write $w$ directly in a power series in $(\xi,b)$.
We will see that $w$ is a sum of terms of
the form $b^m \xi^\nu /|\xi|^{2p}$ with $|\nu| \geq 2p$, and such a term makes a contribution to the error $\eps$ of order
$\del^{m+|\nu|-2p+1}$.  For this purpose we need a qualitative understanding of the inverse of the mapping $\al \mapsto \xi$.
Thinking of polar coordinates in \eqref{defxi}, we do not change the angle but we make a change along each ray depending on the angle.  Thus it is enough to
consider the inverse of the mapping $|\al| \mapsto |\xi|$.  We will do this using the Lagrange Inversion Theorem \cite{ww,krantz}; the theorem is usually
stated for analytic functions, but for $C^N$ functions it can be applied to the Taylor polynomial.

We start by rewriting  \eqref{rsq} as
\beq \label{xi2al2}
  |\xi|^2/|\al|^2 = 1 +  \sum_{|\nu|=4,5} c_\nu\al^\nu/|\al|^2
  + b\sum_{2\leq|\nu|\leq 4} c'_\nu\al^\nu/|\al|^2 
  + O(|(\al,b)|^4) \eeq
Here $O(|(\al,b)|^4)$ means $O(|\al|^4 + b^4)$.  With 
$ u = \al/|\al| = \xi/|\xi|  $,
we can substitute $\al^\nu = u^\nu|\al|^{|\nu|}$ in \eqref{xi2al2}.  We then regard \eqref{xi2al2} as a power series in $|\al|$ in which 
the coefficients depend on $b$ and $u$.  We will say that such a series is of type A if the coefficient of the $n$th power is a
polynomial in $b$ and $u$ with terms $u^\nu$ such that $|\nu| - n$ is even.  Then \eqref{xi2al2} is of type A.  Multiplication of series
preserves type A; thus powers of $|\xi|^2/|\al|^2$ have series of type A.  We note that the $k$th term in a product series depends only on the first $k$ terms in the factors.
Using the power series for $(1+x)^{-1/2}$ we can write a similar expression for $|\al|/|\xi|$ with terms as in \eqref{xi2al2} and their products.  
This series is also of type A;
the same is true for powers of $|\al|/|\xi|$. 
We now apply the Lagrange Theorem to the function $F(|\al|) = |\xi|$.  According to the theorem, $F^{-1}$ has a series in $|\xi|$, with remainder,
 such that the coefficient of $|\xi|^n$ is proportional to the
coefficient of $|\al|^{n-1}$ in the series for 
$\left(|\al|/F(|\al|)\right)^n$ =
$(|\al|/|\xi|)^n$.
This quantity has factors $u^\nu$ with $|\nu| - (n-1)$ even.
We now divide this expression for $F^{-1}(|\xi|) =|\al|$ by $|\xi|$ so that the earlier parity is restored.
We have shown that $|\al|/|\xi|$ has a series in $|\xi|$ which is type A.
Finally we rewrite $u^\nu|\xi|^n$ as $ \xi^\nu|\xi|^{n-|\nu|}$, and
in summary we have shown that
\beq \frac{|\al|}{|\xi|} = 
        \sum c_{m \nu p} b^m \frac{\xi^\nu}{|\xi|^{2p}}
  + O(|(\xi,b)|^4) \eeq
where $ m \geq 0$, $|\nu| \geq 2p$, and
$  m + |\nu| - 2p \leq 3$. With $\al_j = (|\al|/|\xi|)\xi_j$
we get a similar expression for $\al$ as a function of $\xi$.

The function $f(\al)$ and the factor $|T_1\times T_2|$ in $w$ have series in $\al$ which can be converted to $\xi$.   The Jacobian
is
\beq \left| \frac{\pa\al}{\pa\xi} \right| = 
  \mu^2 + \mu\xi\frac{\pa\mu}{\pa|\xi|}\,, \quad
    \mu = \frac{|\al|}{|\xi|}  \eeq
It has terms of the same type as those in $|\al|/|\xi|$.
The Jacobian has leading term $1$ and is bounded but not smooth as
$\xi \to 0$.
We conclude that $w$ has the expression 
\beq \label{wsglexpand}
  w(\xi,b) = \sum c_{m \nu p} b^m \frac{\xi^\nu}{|\xi|^{2p}}
  + R(\xi,b) \eeq
where $m\geq 0$, $|\nu| \geq 2p$, $m + |\nu| - 2p \leq 3$, and 
$R(\xi,b) = O(|(\xi,b)|^4) $.

To find the contribution $\eps_{m\nu p}$ to the error \eqref{epsxib}
from a term in \eqref{wsglexpand} with a particular 
$(m,\nu,p)$ we will integrate in polar coordinates.  The angular integral is zero by symmetry unless $\nu_1$, $\nu_2$ are both even.  Let $n = |\nu| - 2p$, the degree of $\xi$.  With the restriction $m + n \leq 3$ the possible nonzero terms have
$n = 0$ and $0 \leq m \leq 3$ or $n = 2$ with $m = 0,1$.  To carry out the integration, we rescale variables to $\xi = \del\zeta$,
$b = \del\lam$, and write $\zeta$ in polar coordinates.  With
$s = |\zeta|$ we obtain
\beq \eps_{m\nu p} = c_{m\nu p} b^m \del^{n+1} I_n(\lam)
                 = c_{m\nu p} \lam^m \del^{m+ n+1} I_n(\lam) \eeq
where
\beq I_n(\lam) = \int_0^\infty 
    \frac{\erfc(\sqrt{s^2 + \lam^2})}
       {\sqrt{s^2 + \lam^2}} s^{n+1}\,ds  \eeq
In a similar way we see that the remainder $R$ leads to an error which is $O(\del^5)$.  In summary we can express the error as
\beq \label{epsform} 
  \eps =  \del p_0(b)I_0(\lam) + \del^3 p_2(b)I_2(\lam)
  + O(\del^5) \eeq
where $p_0, p_2$ are polynomials in $b$ with $\deg p_0 = 3$, 
$\deg p_2 = 1$.  They depend only on the surface and $b$, not
$\del$ or $\lam$.  For fixed $b$ and $h$ they are unknown coefficients.  To normalize the equation we set $\del = \rho h$
and rewrite it as 
\beq \label{formnorm}
  \eps =  c_1\rho I_0(\lam) + c_2\rho^3 I_2(\lam)
  + O(\del^5) \eeq
This conclusion is equivalent to \eqref{form5}, which we use with three choices of $\del$ to solve for the single layer potential within $O(\del^5)$.

For the double layer potential, in view of \eqref{dblreg} and \eqref{dblsub}, we can write the error from regularizing as
\beq \eps = \frac{1}{4\pi}\int \phi(r/\del)
  \frac{(\xx(\al)-\yy)\cdot \nn(\al)}{r^3}(g(\al)-g(0))\,dS(\al) \eeq
where
\beq \phi(r) = - \erfc(r) - (2/\sqrt{\pi})re^{-r^2} \eeq
and after changing from $\al$ to $\xi$,
\beq \eps = \frac{1}{4\pi}\int
 \frac{\phi(\sqrt{|\xi|^2 + b^2}/\del)}{(|\xi|^2 + b^2)^{3/2} }
    w(\xi,b)\,d\xi  \eeq
where now
\beq  w(\xi,b) = [(\xx-\yy)\cdot \nn][g(\al) - g(0)]
\left|\frac{\pa\al}{\pa\xi}\right| |T_1\times T_2| \eeq
We find
\beq (\xx-\yy)\cdot \nn = -b + O(|(\xi,b)|^2) \eeq
and note $g(\al) - g(0) = O(|\xi|)$.  Thus each term in $w$ now has at least two additional factors.  We expand $w$ as in \eqref{wsglexpand}
but now include terms with $m + n \leq 5$, where again
$n = |\nu| - 2p$.  The term $(m,\nu,p)$ now contributes an
error of order $\del^{m+n-1}$, rather than $\del^{m+n+1}$ as before.
From the last remark, each nonzero term must have $m \geq 1$ and $n \geq 1$ or $m = 0$ and $n \geq 3$. By symmetry a term that contributes a nonzero error must have $m \geq 1$ and $n \geq 2$
or $m = 0$ and $n \geq 4$. The possible terms with $m+n \leq 5$
are $(m,2)$ with $m = 1,2,3$ and $(m,4)$ with $m = 0,1$.
Rescaling the integrals we find
\beq \label{formdbl}
  \eps =  \del p_0(b)J_0(\lam) + \del^3 p_2(b)J_2(\lam)
  + O(\del^5) \eeq
with $\deg p_0 = 3$, $\deg p_2 = 1$, and
\beq J_n(\lambda) = \int_0^\infty
\frac{\phi(\sqrt{s^2+\lam^2})}{(s^2+\lam^2)^{3/2}}s^{n+3}\,ds  \eeq
In fact
\beq  J_0 = - 2I_0\,, \qquad J_2 = -4I_2 \eeq
so that \eqref{formdbl} is equivalent to \eqref{epsform}, and we can solve for the double layer as in \eqref{form5}.

The expansions can be carried further in the same manner.  For the single layer integral we can refine the error expression \eqref{epsform} to
\beq \eps =  \del p_0(b)I_0(\lam) + \del^3 p_2(b)I_2(\lam)
     + \del^5 p_4(b)I_4(\lam) + O(\del^7) \eeq
For the double layer \eqref{formdbl} is replaced by
\beq \eps =  \del p_0(b)J_0(\lam) + \del^3 p_2(b)J_2(\lam)
     + \del^5 p_4(b)J_4(\lam) + O(\del^7) \eeq
Each of these expressions leads to a system of four equations in four unknowns, using four different choices of $\del$.
In fact $J_4 = -6I_4$, so that again we may use the same equations for both cases.

For the Stokes single layer integral, calculated in the form
\eqref{stosglsub}, \eqref{Sreg}, the first term is equivalent to the single layer potential \eqref{sgllayer}.  The second term resembles the double layer \eqref{dbllayer}.  We note the integrand has a factor
$\fftw(\xx)\cdot(\yy - \xx)$ with
$ \fftw(\xx) = \ff(\xx) - (\ff(\xx_0)\cdot\nn(\xx_0))\nn(\xx) $.
Thus $\fftw \cdot\nn = 0$ at $\xx = \xx_0$, and since
$\yy - \xx = b\nn(\xx_0) + O(\xi)$, the numerator of the integrand is $O(\xi)$.  The discussion above for the double layer now applies to this second term, leading to the same expression for the error.

For the Stokes double layer integral, with the subtraction
\eqref{stodblsub} and the kernel rewritten as in \eqref{Treg}, the first term is again like the harmonic double layer.  For the second term, regularized with $s_3$, the numerator in the expansion will have terms $O(\xi^3)$ or higher.  By symmetry the terms that contribute nonzero error have $O(\xi^4)$ or higher.  We get an expansion for the error in the second term in the form
\beq 
  \eps =  \del p_0(b)K_0(\lam) + \del^3 p_2(b)K_2(\lam)
  + O(\del^5) \eeq
with
\beq K_n(\lambda) = \int_0^\infty
\frac{\phi_3(\sqrt{s^2+\lam^2})}{(s^2+\lam^2)^{5/2}}s^{n+5}\,ds  \eeq
and $\phi_3 = 1 - s_3$.
We find that $K_0 = (8/3)I_0$ and $K_2 = 8I_2$, so that once again we can use \eqref{form5} for extrapolation.

\section{
Surface quadrature and the discretization error}

We use a quadrature rule for surface integrals introduced in \cite{wilson-10} and used in \cite{b04,byw,tbjcp}.  We cover the surface with a three-dimensional grid with spacing $h$.  The quadrature points have the form $\xx = (ih,jh,x_3)$, i.e., points on the surface $\Gamma$ whose projections on the $(x_1,x_2)$ plane are grid points, and similarly for the other two directions.  We only use points for which the component of the normal vector in the distinguished direction is no smaller than $\cos\theta$ for a chosen angle $\theta$.  In our case we take $\theta = 70^o$.  The weights are determined by a partition of unity $\psi_1, \psi_2, \psi_3$ on the unit sphere; it is applied to the normal vector at each point.  

\begin{figure}[!htb]
\centering
\scalebox{0.65}{\includegraphics{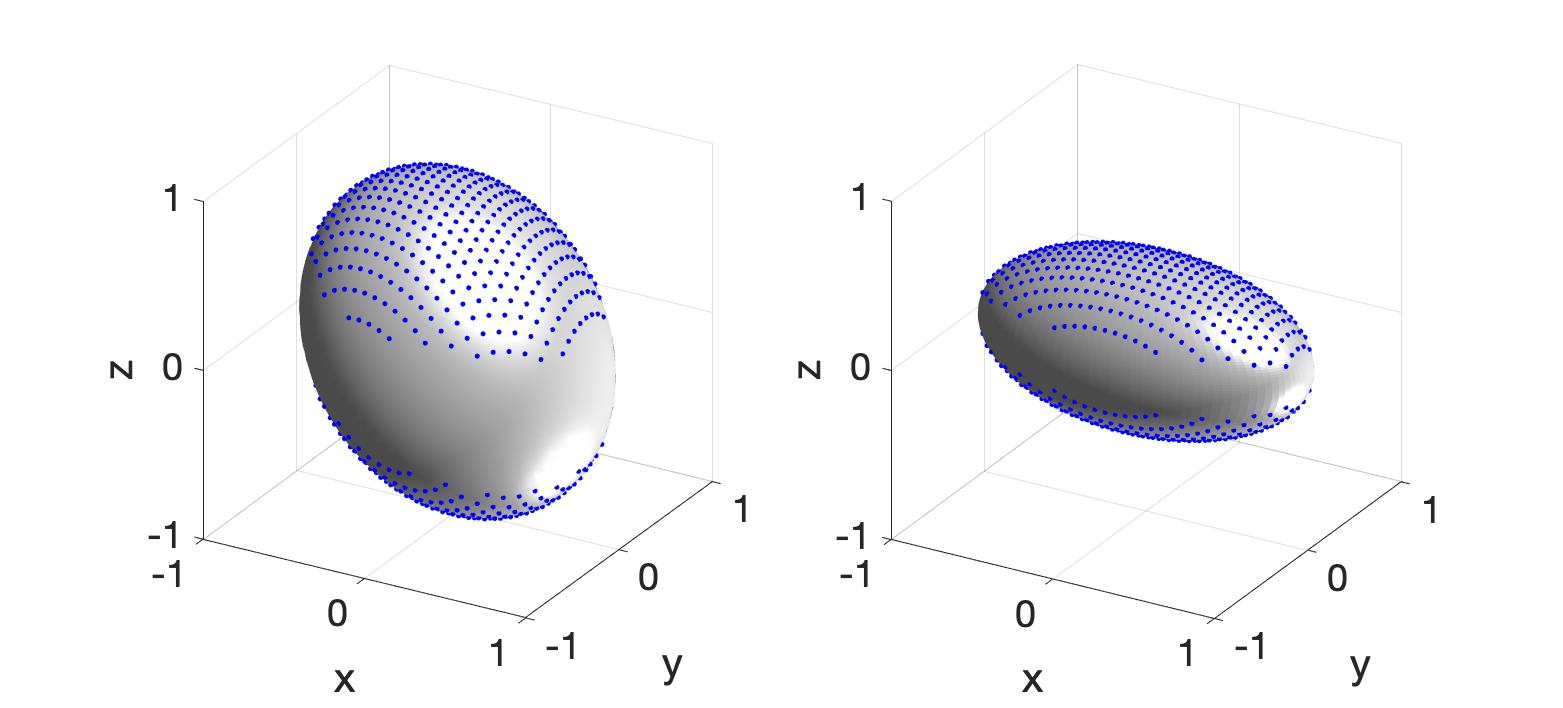}} 
\caption{The rotated (1,.8,.6) ellipsoid (left) and the (1,.5,.5) spheroid (right).}
\label{figure:ellipsoids}
\end{figure}

We define three sets of quadrature points $\Gamma_1, \Gamma_2, \Gamma_3$ as
\beq \Gamma_3 = \{\xx = (ih,jh,x_3) \in \Gamma : |n_3(\xx)| \geq \cos\theta\}  \eeq
where $n_3$ means the third component of the normal vector, and similarly for $\Gamma_1,\Gamma_2$. The quadrature points of the set $\Gamma_3$ are shown for two ellipsoids in Figure~\ref{figure:ellipsoids}.
To construct the partition of unity we start with the bump function
\beq b(r) = \exp(a r^2/(r^2-1))\,, \quad |r|<1\,; \qquad
    b(r) =  0\,, \quad |r| \geq 1              \eeq 
Here $a$ is a parameter.
For a unit vector $\nn = (n_1,n_2,n_3)$
we define 
\beq \beta_i(\nn) = b(\cos^{-1}|n_i|)/\theta) \,, \quad \psi_i(\nn) = \beta_i(\nn)/\left(\sum_{j=1}^3 \beta_j(\nn)\right)  \eeq
The quadrature rule for a surface integral with integrand $f$ is
\beq \int_\Gamma f(\xx)\,dS(\xx) \,\approx\, \sum_{i=1}^3 \sum_{\xx \in \Gamma_i} f(\xx)w_i(\xx) \,,
    \qquad w_i(\xx) = \frac{\psi_i(\nn(\xx))}{|n_i(\xx)|} \,h^2 \eeq
It has high order accuracy as allowed by the smoothness of the surface and the integrand.
The weights cut off the sum in each plane, and each sum has the character of
the trapezoidal rule without boundary; see \cite{wilson-10}.

In earlier work we chose the parameter $a$ to be $1$.  Here we use $a = 2$.  We have found from error estimates in \cite{neglect}, discussed below, as well as numerical experiments, that the discretization error is controlled better with this choice.  We do not recommend using $a > 2$ because of increased derivatives.

The full error in this method consists of the regularization error plus the discretization error; symbolically
\beq \textstyle\sum_\del \,-\, \smint \,=\, \left( \smint_\del \,-\, \smint \right) \,+\, 
  \left(\textstyle\sum_\del \,-\, \smint_\del \right)  \eeq
For either the single layer potential \eqref{sgllayer} or the double layer \eqref{dblsub} the discretization error
arbitrarily close to $\Gamma$
 can be written as 
\beq \label{derr} 
  c_1 h \,+\, C_2 h^2 \exp(-c_0 \del^2/h^2) \,+\, O(\del^5) \,,
\quad c_1 = c_1(b/\del,\del/h) \eeq
which at first appears inaccurate.  Formulas for the first term were given in \cite{b04,byw}, based on approximating the surface locally as a plane.  They can be used as corrections.  Estimates for these formulas were given in \cite{neglect}.  With
the parameter choices here, in particular with $\del/h \geq 2$, it was shown that
\beq c_1^{(S)} \leq 2.1\cdot 10^{-7}\max{|f|}\,, \quad
   c_1^{(D)} \leq 8.3\cdot 10^{-7}\max{|\nabla g|}  \eeq
for the single and double layer respectively, and they decrease rapidly as $\del/h$ increases.  Here $\nabla g$ means the tangential gradient.
The $h^2$ term in \eqref{derr} evidently decreases rapidly as $\del/h$ increases, as does $c_1$.  With $\theta = 70^o$, $c_0 \approx 1.15$;
see \cite{byw}, Sect. 3.4.  However $C_2$ depends on the surface and integrand and could be large.  With moderate resolution we expect that the
discretization error is controlled by the regularization.  If desired the formulas for $c_1h$ in \cite{byw}
could be used as corrections with the present method; they are infinite series, but only the first few terms are significant.  To ensure that the regularization error dominates the discretization error for small $h$ we can choose $\del$ proportional to $h^q$, with $q < 1$, so that
$\del/h$ increases as $h \to 0$.


\section{Numerical examples}
We present examples computing single and double layer integrals at grid points within $O(h)$ of a surface, for harmonic potentials and for Stokes flow.  The points are selected from the three-dimensional grid with spacing $h$ which determines the quadrature points on the surface, as described in Sect. 4.  With the fifth order regularization the results are in general agreement with the theoretical predictions.  With moderate resolution and
$\delta/h$ constant the errors are about $O(h^5)$.  With $\delta$ proportional to $h^{4/5}$ the error is about $O(h^4)$.  For the harmonic potentials we also test the seventh order method; the errors are typically smaller but the order of accuracy is less predictable. 
It is likely that the discretization error is relatively more significant with the smaller errors of the seventh order case.

We report maximum errors and $L^2$ errors, defined as
\beq \|e\|_{L^2} = \left( \sum_{\yy} |e(y)|^2 / N \right)^{1/2} \eeq
where $e(y)$ is the error at $\yy$ and $N$ is the number of points. 
We present absolute errors; for comparison we give approximate norms of the exact solution.

{\bf Harmonic Potentials.}
We begin with known solutions on the unit sphere.  We test the single and double layer separately.  We compute the integrals at grid points first within distance $h$ and then on shells at increasing distance.  In the latter case we also find values computed without regularization.
We then compute known harmonic functions on three other surfaces which combine single and double layers.  

The single and double layer potentials, \eqref{sgllayer} and \eqref{dbllayer}, are harmonic inside and outside the surface~$\Gamma$. They are characterized by the jump conditions 
\beq [\Sl(\xx)] = 0\,,\quad [\pa \Sl(\xx)/\pa\nn] = f(\xx) \eeq
\beq [\Dl(\xx)] = -g(\xx)\,,\quad [\pa \Dl(\xx)/\pa\nn] = 0 \eeq
where $[\cdot]$ means the value outside $\Gamma$ minus the value inside.

For the unit sphere we use the spherical harmonic function
\beq f(\xx) = 1.75(x_1 - 2x_2)(7.5x_3^2 - 1.5)\,, \quad |\xx| = 1 \eeq
for both the single and double layer integrals.  The functions
\beq u_-(\yy) = r^3f(\yy/r)\,,\quad u_+(\yy) = r^{-4}f(\yy/r)\,,
    \quad\; r = |\yy| \eeq
are both harmonic.  We define
$\Sl(\yy)$ by \eqref{sgllayer} and $\Dl(\yy)$ by \eqref{dbllayer}
with $g = f$.  They are determined by the jump conditions,
\beq \Sl(\yy) = - (1/7)u_-(\yy)\,, \quad |\yy|<1\,; \quad\;
   \Sl(\yy) = - (1/7)u_+(\yy)\,, \quad |\yy|>1\eeq
\beq \Dl(\yy) =  (4/7)u_-(\yy)\,, \quad |\yy|<1\,; \quad\;
   \Dl(\yy) = - (3/7)u_+(\yy)\,, \quad |\yy|>1\eeq

\begin{figure}[!htb]
\centering
\scalebox{0.5}{\includegraphics{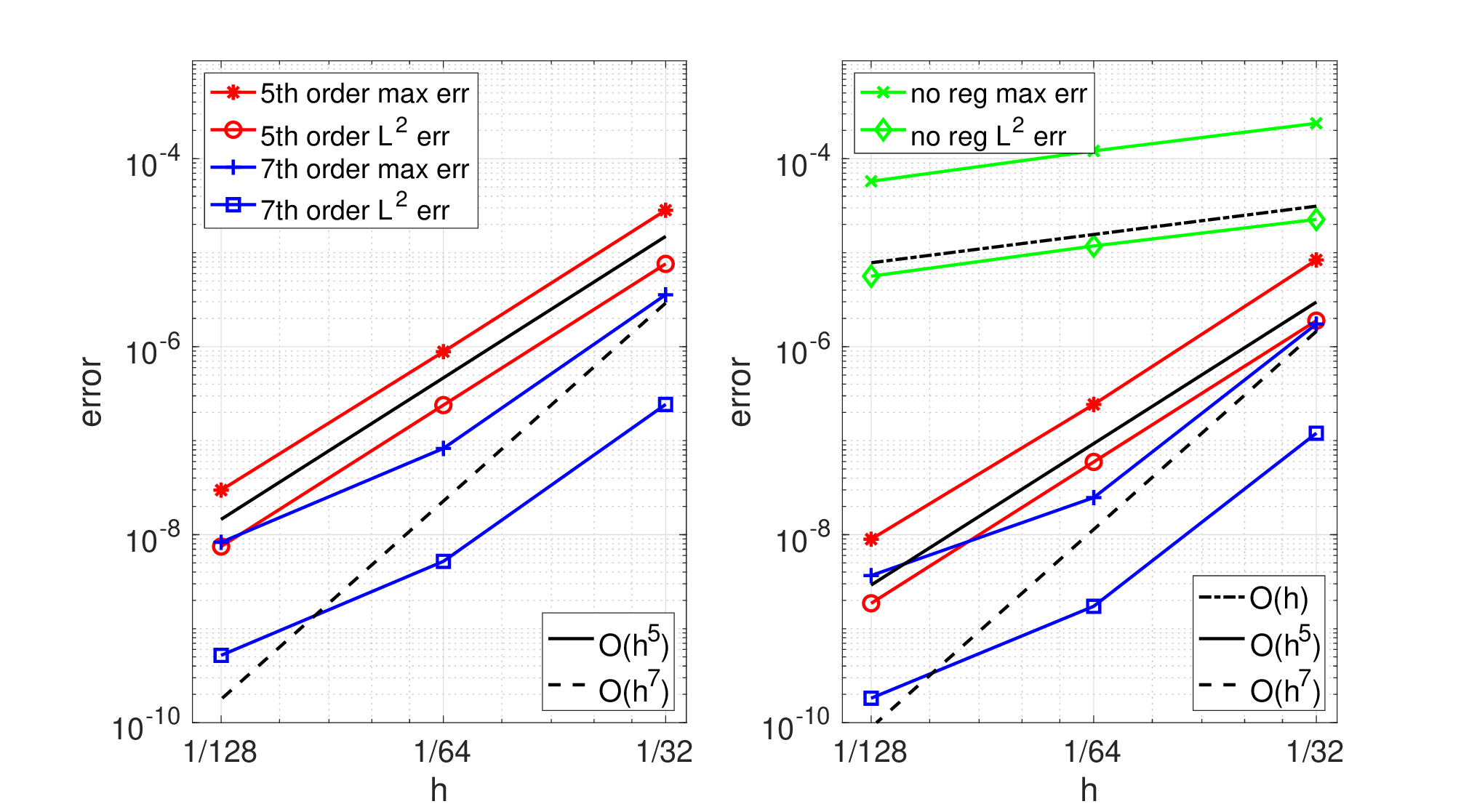}} 
\caption{Errors for the single layer potential on the unit sphere, 
(left) at grid points within distance $h$, computed with the 5th and 7th order regularization, and 
(right) evaluated at distance between $h$ and $2h$, without regularization and with the 5th and 7th order methods.}
\label{figure:tables1_2}
\end{figure}

\begin{figure}[!htb]
\centering
\scalebox{0.5}{\includegraphics{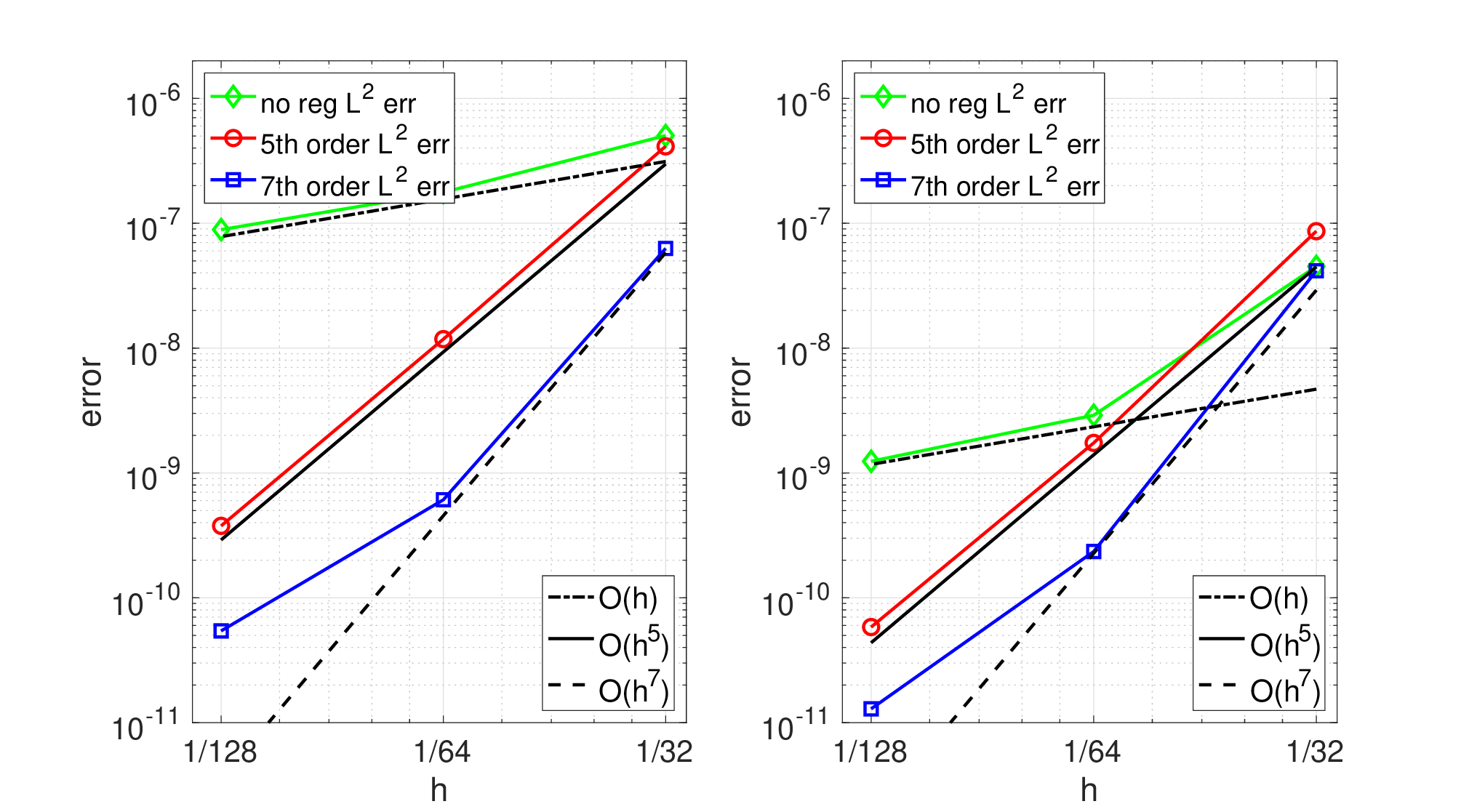}} 
\caption{$L^2$ errors in the single layer potential on the unit sphere, evaluated at distance between $2h$ and $3h$ (left) or $3h$ and $4h$ (right).}
\label{figure:table3}
\end{figure}

We present errors for the single and double layer potentials at grid points at various distances from the sphere.  We begin with the single layer.  We compute the integral as in \eqref{sglreg} and extrapolate as in \eqref{form5}.  Near the sphere the maximum of $|\Sl|$ is about $1.15$ and the $L^2$ norm is about $.50$.  Figure~\ref{figure:tables1_2}, left, shows the $L^2$ and maximum errors for grid points within distance $h$ of the sphere, using fifth or seventh order extrapolation.  For the fifth order we take $\delta/h = 2,3,4$ as previously described, and for
the seventh order we take $\delta/h = 2,3,4,5$. The expected order of accuracy is evident in the fifth order case; the seventh order method has somewhat smaller errors but does not have a discernible order of accuracy, probably because the discretization error is significant.  In subsequent figures we display the errors at nearby grid points at distance between $mh$ and $(m+1)h$ from the sphere, both inside and outside, for $m = 1,2,3$.  We compute the integral with no regularization as well as the fifth and seventh order methods.  Figure~\ref{figure:tables1_2}, right, shows errors for $m=1$ and Figure~\ref{figure:table3} for $m = 2$ and $3$.  The values without regularization in Figure~\ref{figure:tables1_2} appear to be about $O(h)$ accurate.  The fifth order method again has the expected order of accuracy at least for $m = 1$ but becomes less steady with distance.  The errors become smaller overall as the distance increases.  Beyond $4h$ the error without regularization is quite small, suggesting that we can discontinue the regularization for points at least $4h$ from the surface.

\begin{figure}[!htb]
\centering
\scalebox{0.5}{\includegraphics{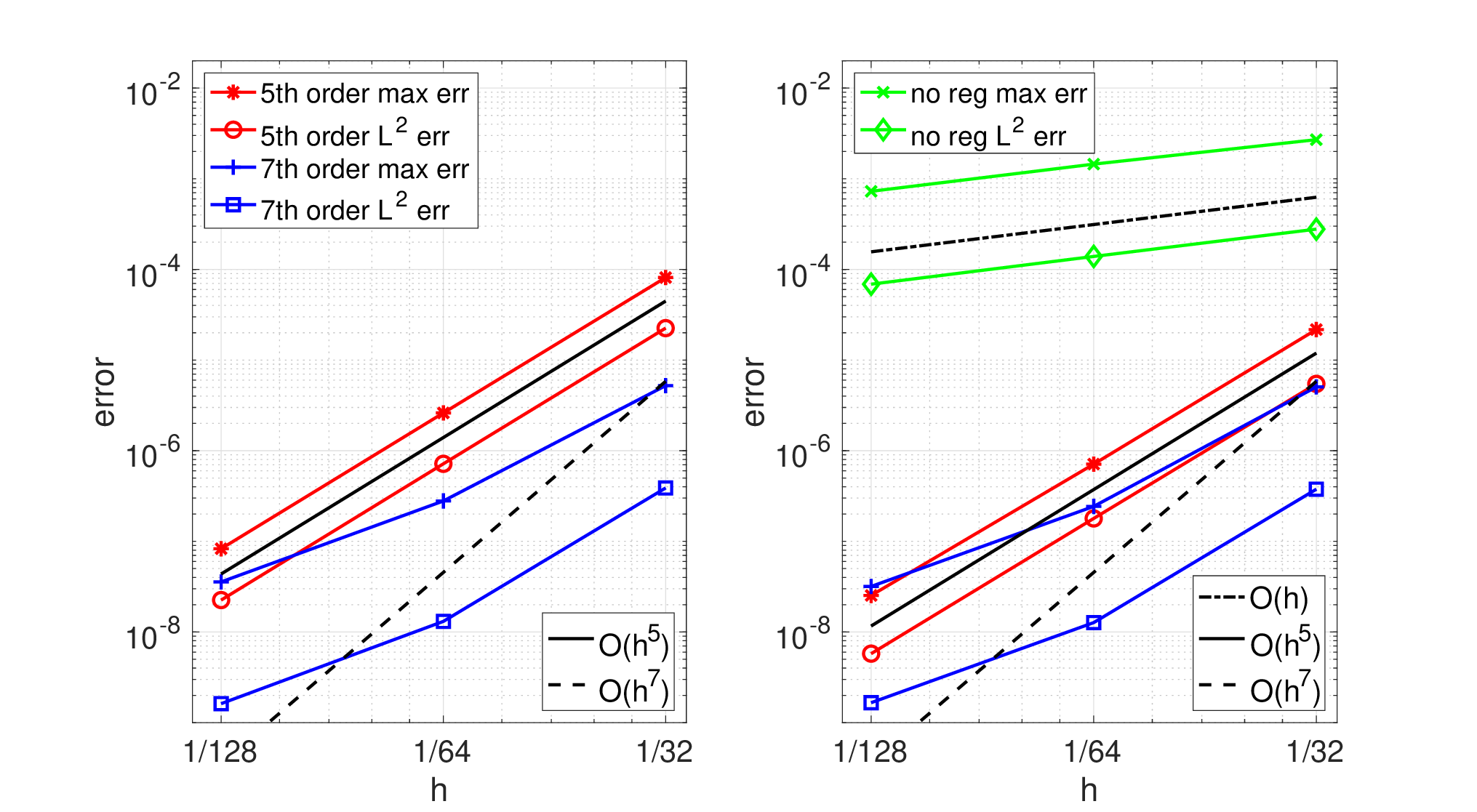}} 
\caption{Errors for the double layer potential on the unit sphere,
(left) at grid points within distance $h$, computed with the 5th and 7th order regularization, and 
(right) evaluated at distance between $h$ and $2h$, without regularization and with the 5th and 7th order methods.}
\label{figure:tables4_5}
\end{figure}

\begin{figure}[!htb]
\centering
\scalebox{0.5}{\includegraphics{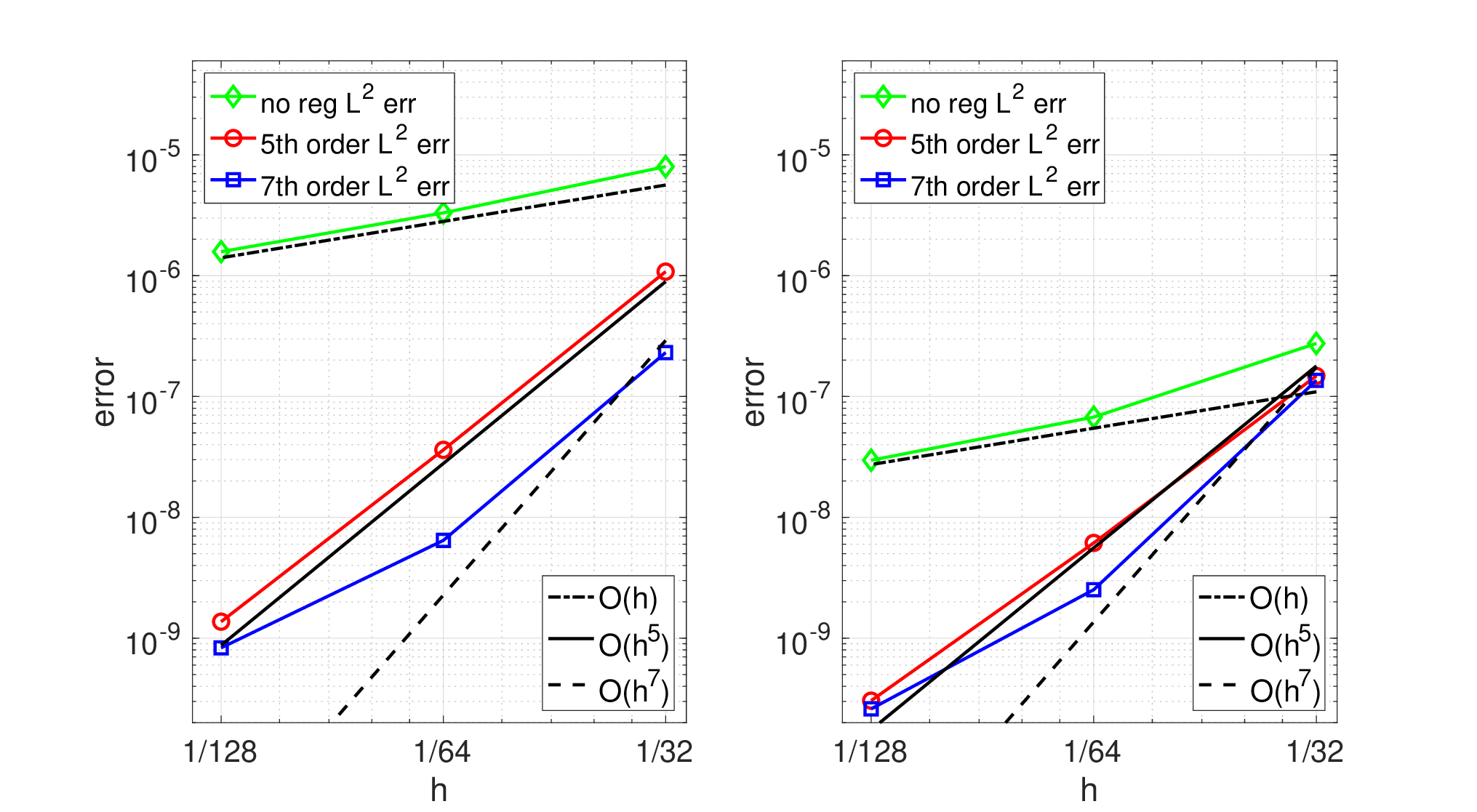}} 
\caption{$L^2$ errors in the double layer potential on the unit sphere, evaluated at
distance between $2h$ and $3h$ (left) or $3h$ and $4h$ (right).}
\label{figure:table6}
\end{figure}

In Figures~\ref{figure:tables4_5},\ref{figure:table6} we present results of the same type for the double layer potential, computed as in \eqref{dblreg}.  They are similar in behavior to those for the single layer.  The maximum of $|\Dl|$ is about $4.6$ and
$\|\Dl\|_{L^2} \approx 1.8$.

For the remaining tests on other surfaces we use a procedure as in \cite{byw} which allows us to have known solutions with an arbitrary surface $\Gamma$.  This provides a test of the single and double layer combined, rather than separately.  We choose harmonic functions $u_+$ outside and $u_-$ inside.  We set $f = -[u]$ and $g = [\pa u/\pa n]$, the jumps across $\Gamma$ as above.  Then assuming $u_+$ decays at infinity, $ u(\yy) = \Sl(\yy) + \Dl(\yy)$ on both sides, where $\Sl$ and $\Dl$ are defined in \eqref{sgllayer}, \eqref{dbllayer}.  We choose
\beq u_-(\yy) = (\sin{y_1} + \sin{y_2})\exp{y_3}\,,
   \quad u_+(\yy) = 0  \eeq

In these tests we again use $\del/h = 2,3,4$ with the fifth order method and $\del/h = 2,3,4,5$ with seventh order.  We also
choose $\delta$ proportional to
$h^{4/5}$ with the fifth order method and $h^{4/7}$ with the seventh order method, so that the predicted order of error is 
$O(h^4)$.  We choose constants so that $\del$ agrees with the earlier choice at $1/h = 64$.

\begin{figure}[!htb]
\centering
\scalebox{0.5}{\includegraphics{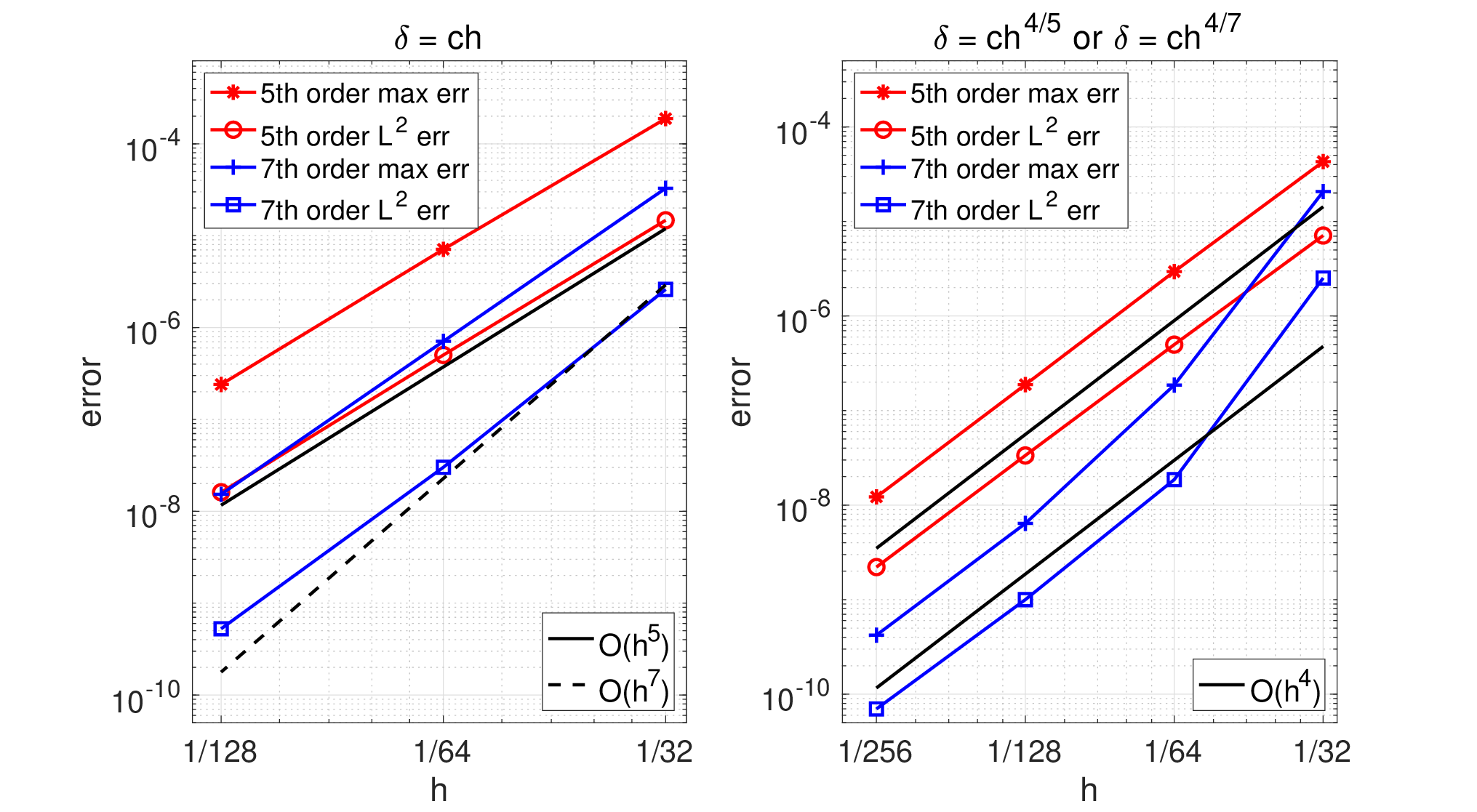}} 
\caption{(Left) Errors for the single and double layers on a rotated ellipsoid at grid points within distance $h$, with the 5th order and 7th order methods, $\del$ proportional to $h$. (Right) Errors for the rotated ellipsoid, at  grid points within distance $h$ in the first octant; $5$th and $7$th order methods with $\del$ chosen to correspond to $O(h^4)$ accuracy.}
\label{figure:tables7_8}
\end{figure}

Our first surface with this procedure is a rotated ellipsoid shown in Figure~\ref{figure:ellipsoids}, left,
\beq \frac{z_1^2}{a^2} + \frac{z_2^2}{b^2} + \frac{z_3^2}{c^2} = 1 \eeq
where $a = 1$, $b = .8$, $c = .6$ and $\zz = M\xx$, where
$M$ is the orthogonal matrix 
\beq M = (1/\sqrt{6})\,[\sqrt{2}\quad 0\quad -2;\;
          \sqrt{2}\quad \sqrt{3}\quad 1;\;
       \sqrt{2}\quad -\sqrt{3}\quad 1].  \eeq
We present results in Figure~\ref{figure:tables7_8}.  In Figure~\ref{figure:tables7_8}, left, we evaluate at
all grid points within distance $h$ with both regularizations.
Figure~\ref{figure:tables7_8}, right, has values at points $y$ within distance $h$
in the first octant, i.e., those with $y_1, y_2, y_3 \geq 0$.
  The accuracy of the fifth order version is close to the prediction; the seventh order version has smaller errors in Figure~\ref{figure:tables7_8}, right, and perhaps approximates the predicted order $O(h^4)$ but not clearly so.
For the left figure the $L^2$ norm of the exact solution is about
$.5$ and the maximum about 1.7.  For the right figure, within the first octant, they are about .76 and 1.4.

\begin{figure}[!htb]
\centering
\scalebox{0.65}{\includegraphics{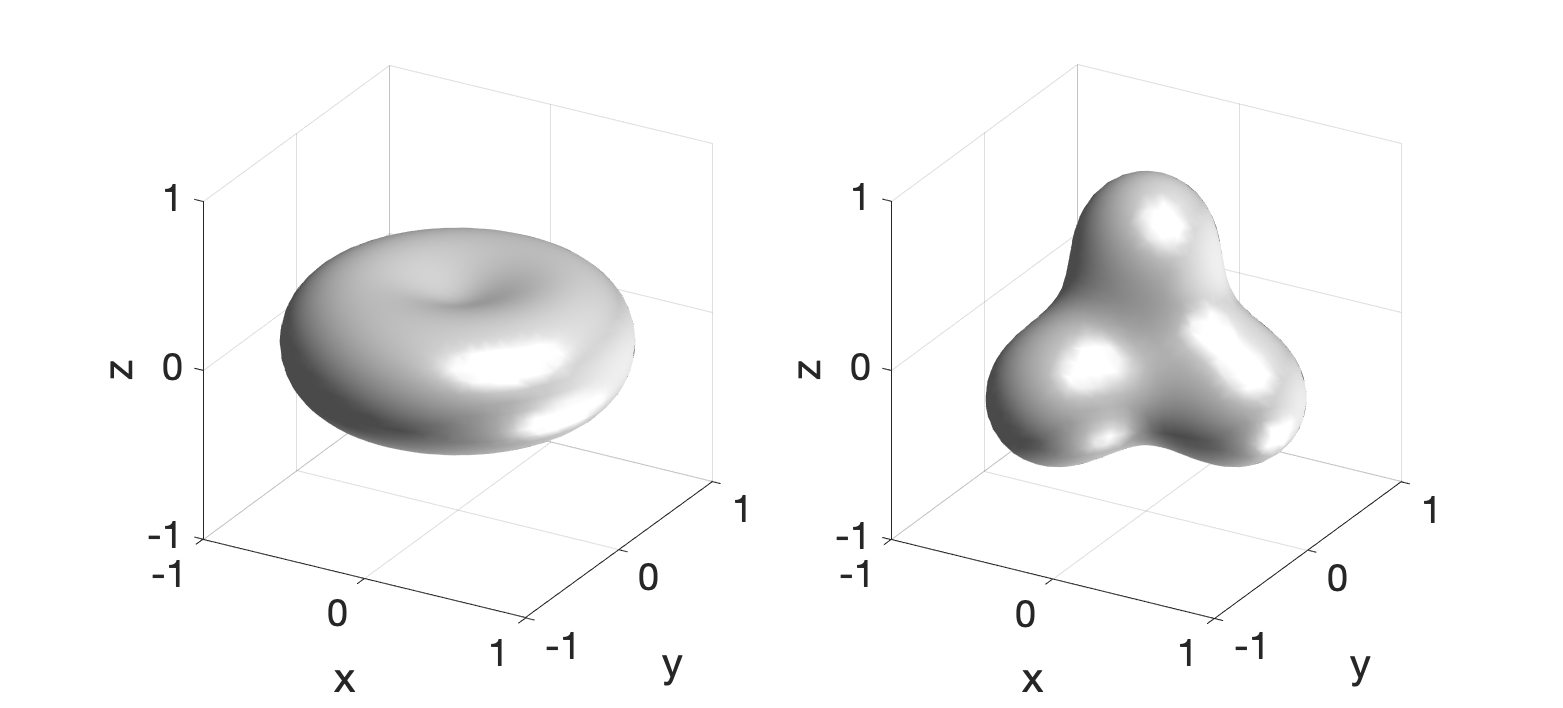}} 
\caption{The Cassini oval surface and the four-atom molecular surface.}
\label{figure:cassini_molecule}
\end{figure}

\begin{figure}[!htb]
\centering
\scalebox{0.5}{\includegraphics{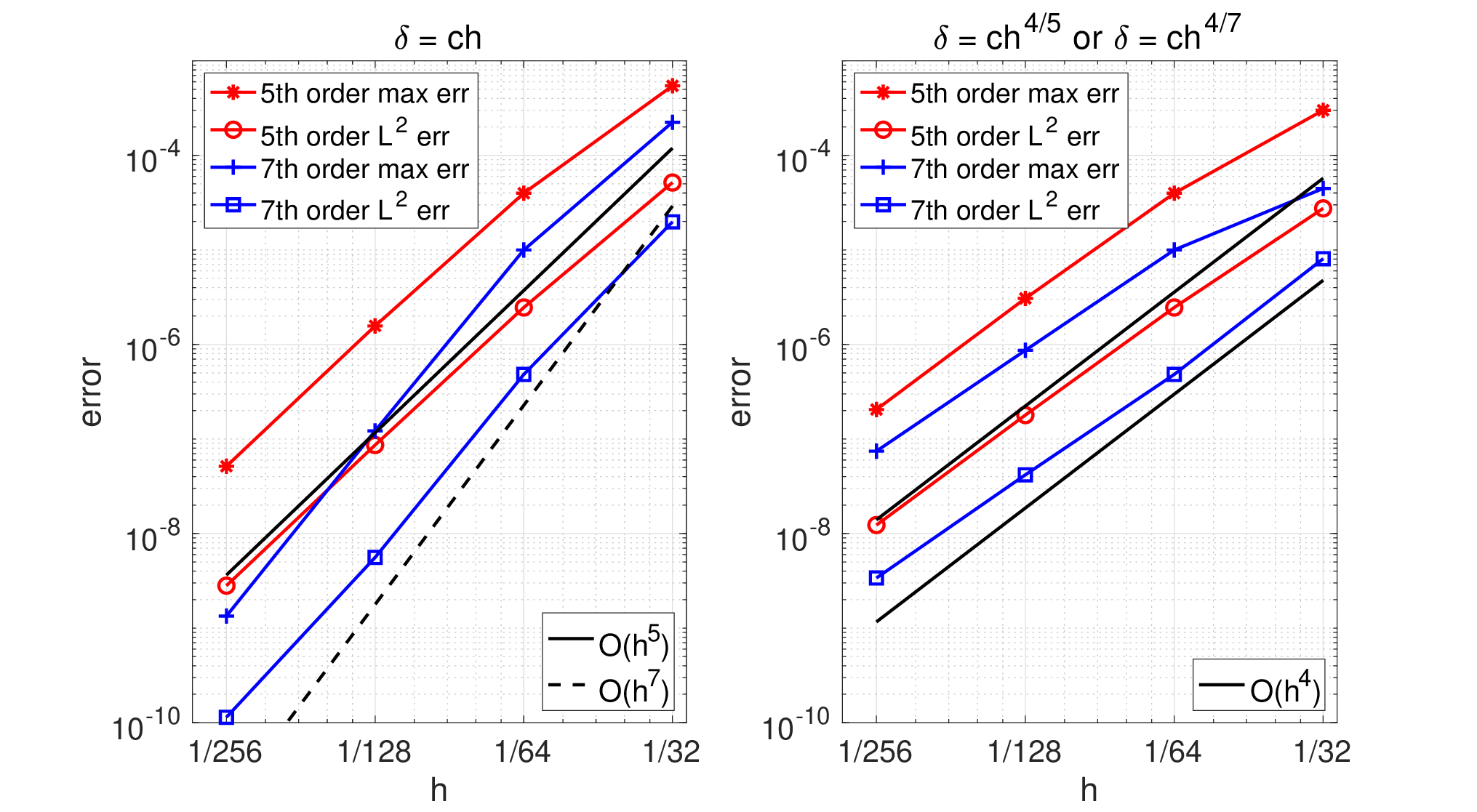}} 
\caption{Errors for the Cassini oval surface, at grid points within distance $h$ in the first octant; $5$th and $7$th order method with $\del$ proportional to $h$ or corresponding to $O(h^4)$ accuracy.}
\label{figure:table9}
\end{figure}

\begin{figure}[!htb]
\centering
\scalebox{0.5}{\includegraphics{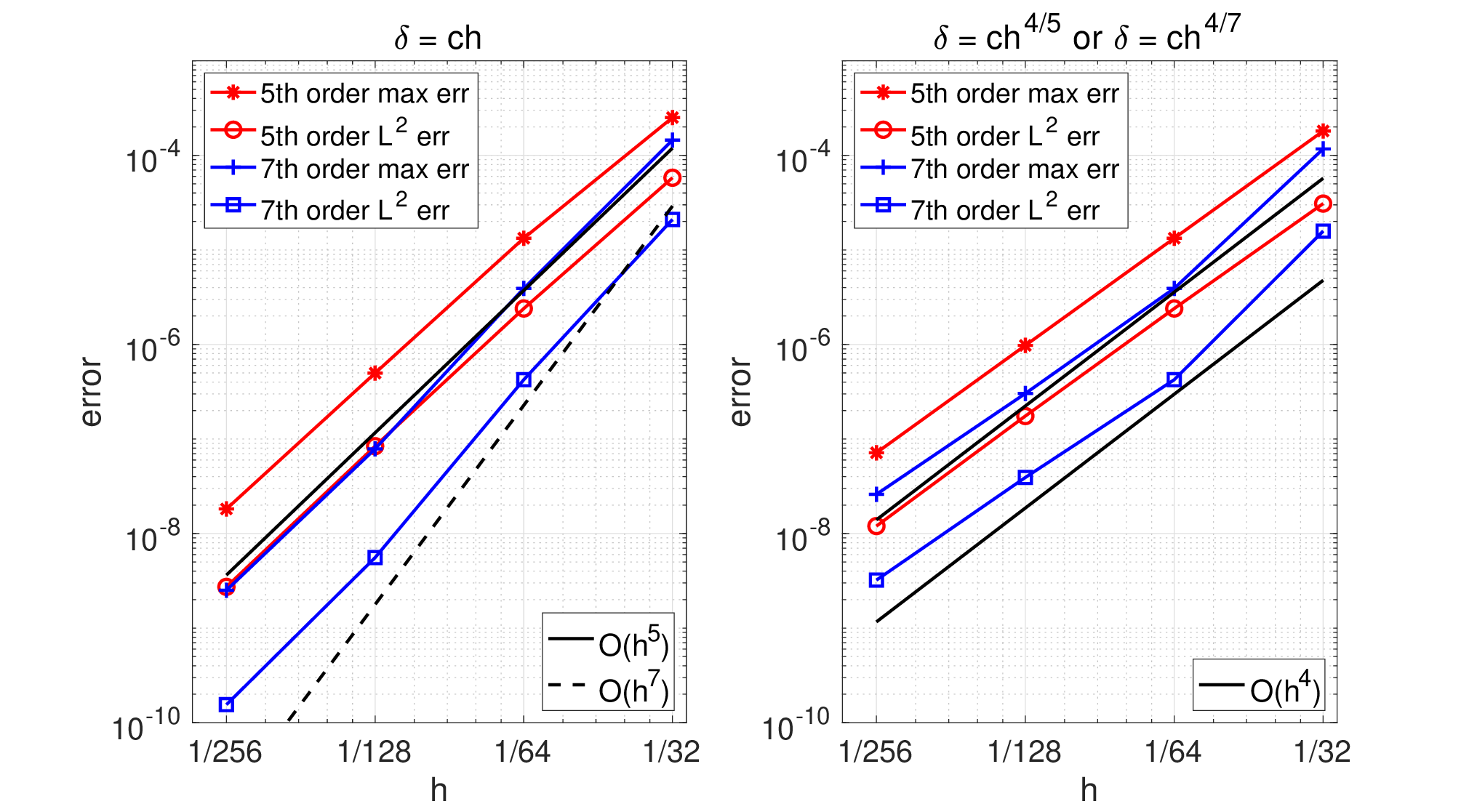}} 
\caption{Errors for the molecular surface, at grid points within distance $h$ in the first octant; $5$th and $7$th order method with $\del$ proportional to $h$ or corresponding to $O(h^4)$ accuracy.}
\label{figure:table10}
\end{figure}

The next example is a surface obtained by revolving a Cassini oval about the $x_3$ axis,
\beq (x_1^2 + x_2^2 + x_3^2 + a^2)^2 - 4a^2(x_1^2 + x_2^2) = b^4 \eeq
with $a = .65$ and $b = .7$.  The final surface represents a molecule with four atoms, 
\beq \label{molesurf} \sum_{i=1}^4 \exp(-|\xx - \xx_k|^2/r^2) = c \eeq
with $r = .5$, $c = .6$, and $\xx_k$ given by
\beq (\sqrt{3}/3,0,-\sqrt{6}/12)\,,\; 
   (-\sqrt{3}/6,\pm .5,-\sqrt{6}/12)\,,\;
   (0,0,\sqrt{6}/4)  \eeq
These surfaces are shown in Figure~\ref{figure:cassini_molecule}.
We compute the solution for grid points in the first octant as before for the ellipsoid, with $\del$ related to $h$ in the same way.  We present errors with fifth or seventh order regularization, with $\delta$ proportional to $h$ or fractional.
The results, reported in Figures~\ref{figure:table9} and~\ref{figure:table10}, are generally similar to those for the rotated ellipsoid.  For both surfaces we see roughly the predicted orders of accuracy in the fifth order case. For seventh order the errors are smaller, but the accuracy in the fractional case is somewhat less than fourth order in $h$.
  For the Cassini surface the $L^2$ norm for the exact values is about $.78$ and the maximum is about $1.45$.  For the molecular surface they are about $.57$ and $1.0$.

\newpage

{\bf Stokes Flow.}  We present examples of three types.  First we calculate the velocity near a translating spheroid in Stokes flow, given as a single layer integral.  We then compute a standard identity for the double layer integral.  Finally we compute a velocity that combines single and double layer integrals on an arbitrary surface, as in the examples above with harmonic potentials.  We have increased $\rho$ to $(3,4,5)$ to make the order of accuracy more evident, even though errors are typically smaller with $(2,3,4)$.  In each case we report errors at grid points within distance $h$ of the surface.

\begin{figure}[!htb]
\centering
\scalebox{0.425}{\includegraphics{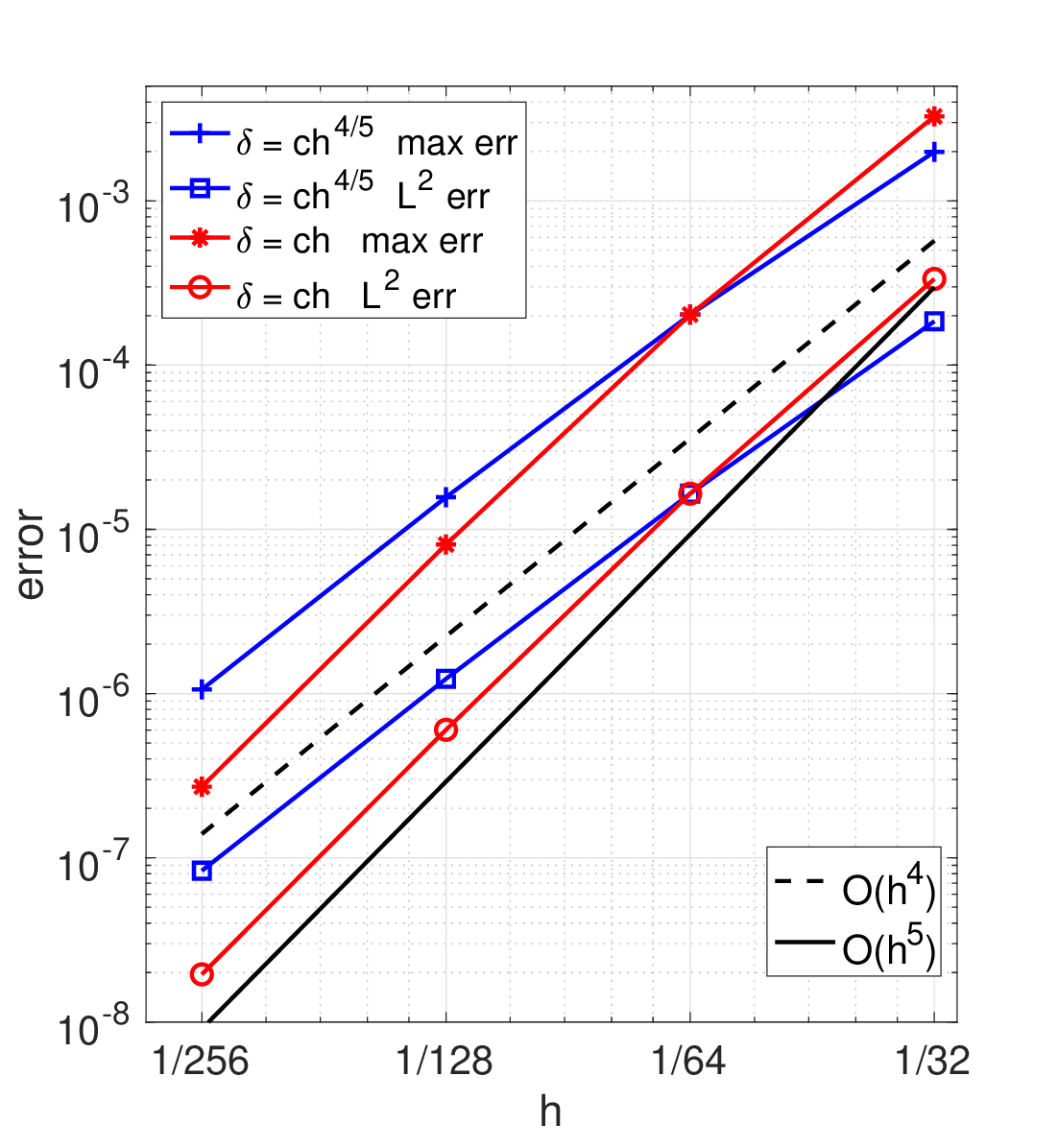}} 
\caption{Errors for the Stokes single layer on a prolate spheroid, at grid points within distance $h$ outside the spheroid.}
\label{figure:table11}
\end{figure}

In our first example we compare the single layer or Stokeslet integral with an exact solution.  We compute the Stokes flow around a prolate spheroid
\beq \label{spheroid}
  x_1^2 + 4x_2^2 + 4x_3^2 = 1 \eeq
 with semi-axes $1,.5,.5$, shown in Figure~\ref{figure:ellipsoids}, right, and translating with velocity $(1,0,0)$.  The fluid velocity is determined by the integral \eqref{SingleLayer}
from the surface traction $\ff$.  Formulas for the solution are given in
\cite{chwang,liron,tbjcp}.  The surface traction is
$$ \ff(\xx) = (f_1(\xx),0,0)\,,\quad f_1(\xx) = 
    \frac{F_0}{\sqrt{1 - 3x_1^2/4}} $$
where $F_0$ is a constant.
We compute the fluid velocity $\uu$ as in \eqref{stosglsub},\eqref{Sreg} and extrapolate as before.  Results are presented in Figure~\ref{figure:table11}.  The exact solution has maximum amplitude~$1$ and $L^2$ norm about~$1$.

\begin{figure}[!htb]
\centering
\scalebox{0.425}{\includegraphics{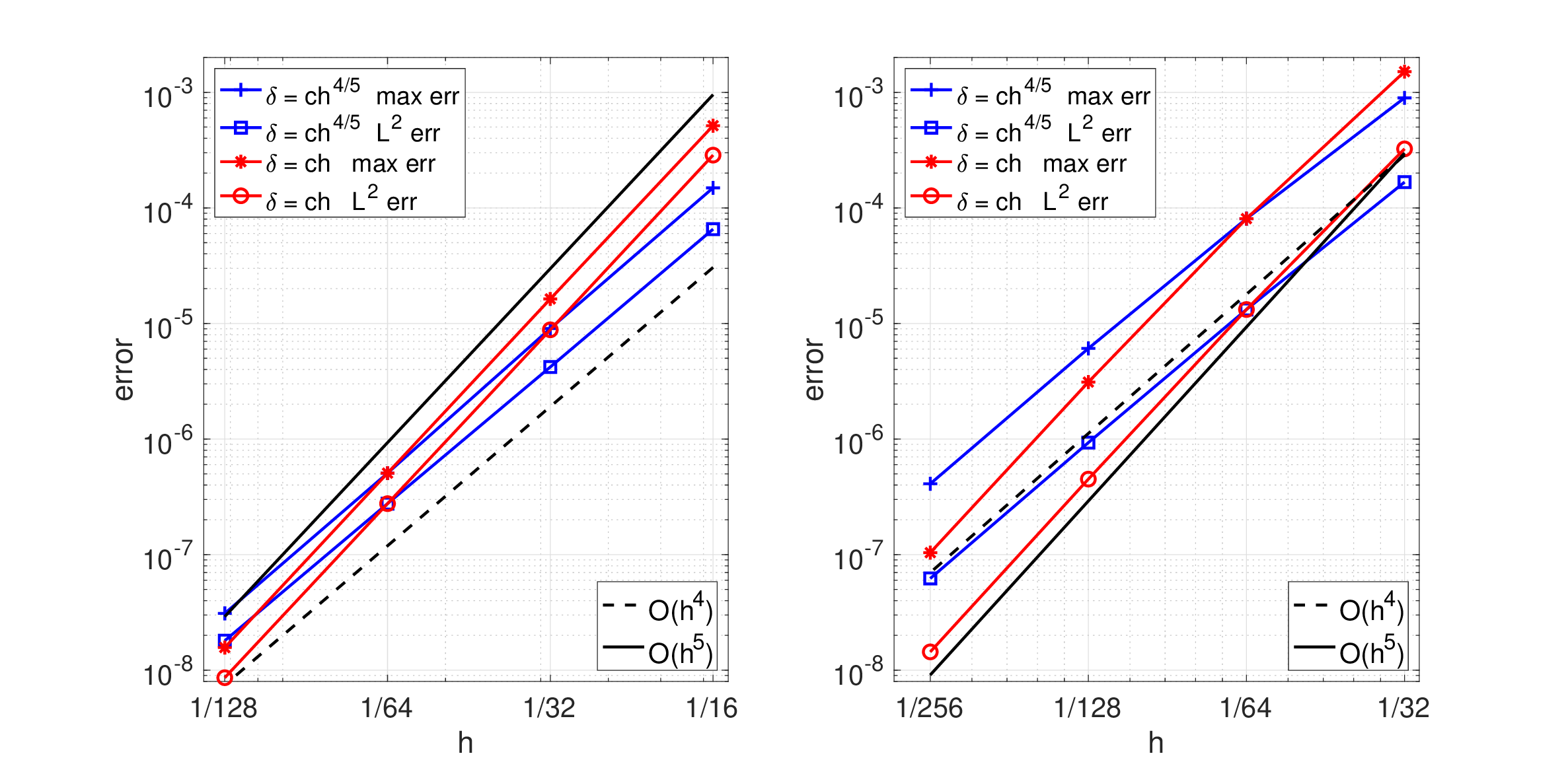}} 
\caption{(Left) Error for the Stokes double layer on the unit sphere, at grid points within distance $h$ on either side of the sphere. (Right) Errors for the Stokes double layer on a prolate spheroid, at grid points within distance $h$ on either side of the spheroid.}
\label{figure:tables12_13}
\end{figure}

Next we test the double layer integral \eqref{DoubleLayer} using the identity (2.3.19) from \cite{pozbook}
\beq 
	\label{DL_identity}
	\frac{1}{8\pi} \epsilon_{jlm} \int_\Gamma x_m T_{ijk} (\bd{x_0,x}) n_k(\bd{x})dS(\bd{x}) = \chi (\bd{x}_0) \epsilon_{ilm} x_{0,m}
\eeq
where $\chi$ = 1, 1/2, 0 when $\bd{x}_0$ is inside, on, and outside the boundary. We set $l=1$ and define $q_j(\bd{x}) = \epsilon_{j1m}x_m = (0,-x_3,x_2)$.  We compute the integral according to \eqref{stodblsub}, \eqref{Tsplit}, \eqref{Treg} and extrapolate.  We report errors for a sphere and for the spheroid \eqref{spheroid} in Figure~\ref{figure:tables12_13}.
For the sphere the maximum value is $1$ and the $L^2$ norm is about $.57$.  For the spheroid the maximum is $ \approx .5$ and the $L^2$ norm is $\approx .3$.

\begin{figure}[!htb]
\centering
\scalebox{0.45}{\includegraphics{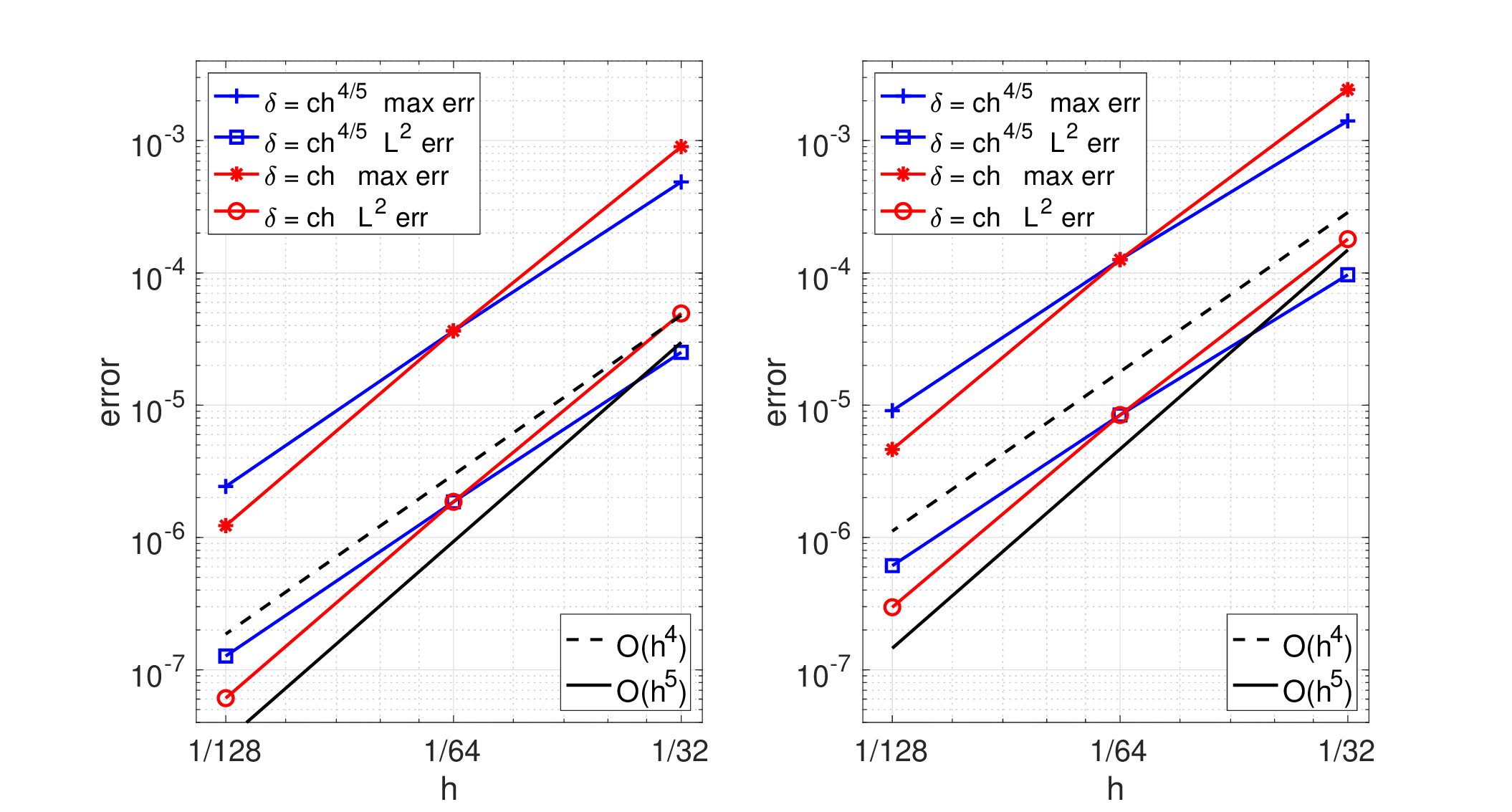}} 
\caption{(Left) Errors for the Stokes single and double layers on the unit sphere, at grid points within distance $h$ on either side of the sphere. (Right) Errors for the Stokes single and double layers on an ellipsoid, at grid points within distance $h$ on either side of the ellipsoid.}
\label{figure:tables14_15}
\end{figure}

\begin{figure}[!htb]
\centering
\scalebox{0.425}{\includegraphics{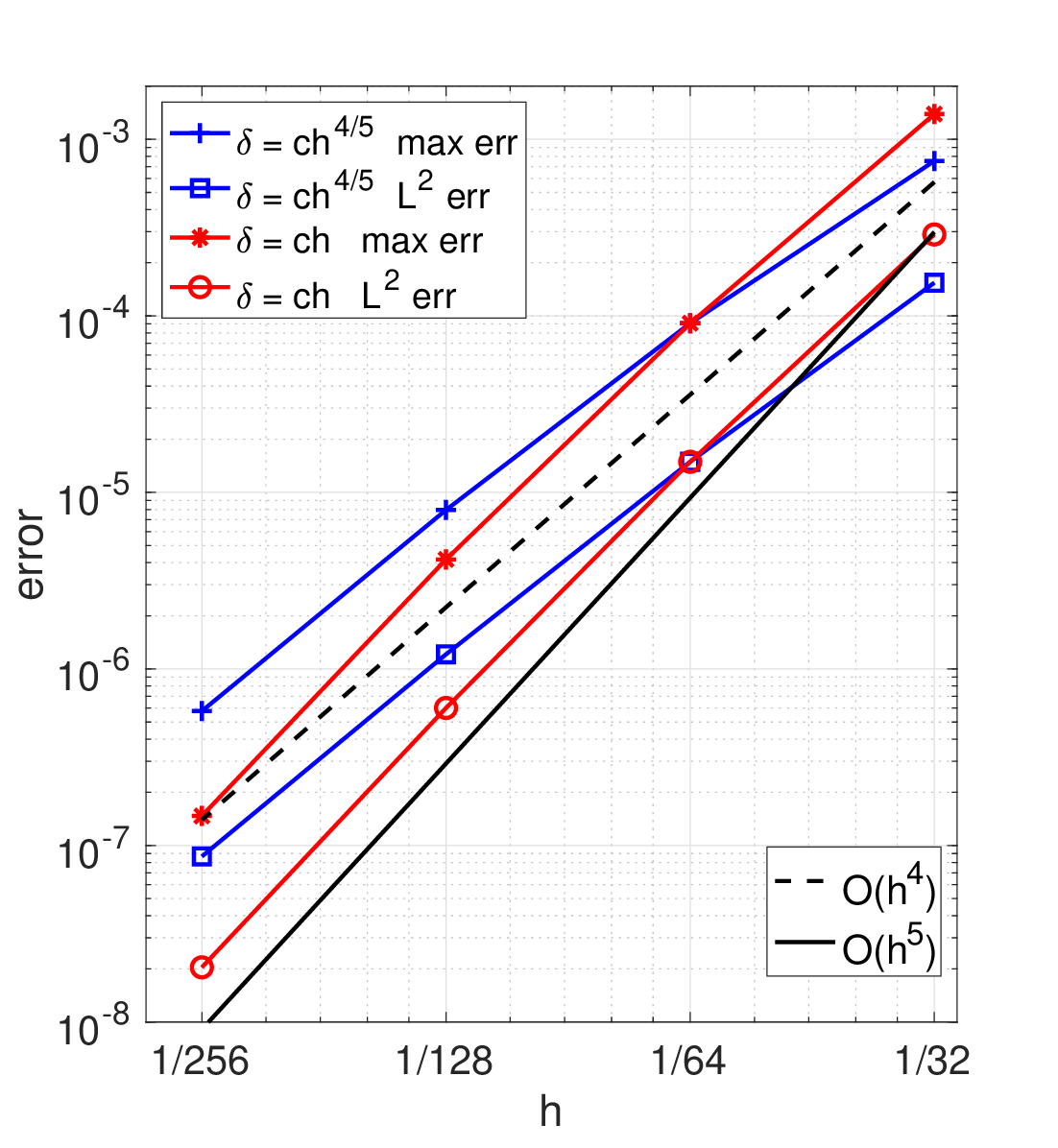}} 
\caption{Errors for the Stokes single and double layers on the four-atom molecular surface, at grid points in the first octant within distance $h$ on either side of the molecule.}
\label{figure:table16}
\end{figure}

In order to test integrals on general surfaces we again use a formula combining the single and double layer integrals.  If $\uu$ is the velocity of Stokes flow outside and inside a surface $\Gamma$, with suitable decay at infinity, then
\beq 
	\label{Sum_SLDL}
	u_i(\bd{y}) = 
	-\frac{1}{8\pi}\int_\Gamma S_{ij} (\bd{y,x}) [f]_j(\bd{x}) dS(\bd{x}) 
	- \frac{1}{8\pi}\int_\Gamma T_{ijk} (\bd{y,x}) [u]_j(\bd{x}) n_k(\bd{x})dS(\bd{x})
\eeq
Here $[f] = f^+ - f^-$ is the jump in surface force, outside minus inside, and $[u]$ is the jump in velocity.  The surface force is the normal stress,
$f^\pm = \sigma^\pm \cdot\nn$, where $\nn$ the outward normal.  The jump conditions are derived e.g. in \cite{pozbook}.  As a test problem we take the inside velocity to be the Stokeslet due to a
point force singularity of strength $\bd{b} = (4\pi,0,0)$, placed at
$\bd{y}_0 = (2,0,0)$. The velocity is
\beq 
	\label{Stokeslet_point}
	u^-_i(\bd{y}) = \frac{1}{8\pi} S_{ij}b_j = \frac{1}{8\pi} \Big( \frac{\delta_{ij}}{r} + \frac{\hat{y}_i \hat{y}_j}{r^3} \Big)b_j
\eeq 
and the stress tensor is 
\beq
	\label{Stress_point}
	\sigma^-_{ik}(\bd{y}) = \frac{1}{8\pi} T_{ijk} b_j = \frac{-6}{8\pi}\frac{\hat{y}_i \hat{y}_j \hat{y}_k}{r^5} b_j
\eeq
where $\hat{\bd{y}} = \bd{y}-\bd{y}_0$, $r=|\hat{\bd{y}}|$. We choose the outside velocity and stress to be zero.  We compute the two integrals in the same manner as above.  We present results for three surfaces: the unit sphere, Figure~\ref{figure:tables14_15}, left; an ellipsoid with semi-axes $1,.8,.6$, Figure~\ref{figure:tables14_15}, right; and the molecular surface \eqref{molesurf}, Figure~\ref{figure:table16}.  For the first two surfaces, the errors are at all grid points within $h$, but for the molecular surface the points are in the first octant only. For the sphere or ellipsoid the maximum velocity magnitude is $\approx 1$ and the $L^2$ norms are $\approx .35$ and $.37$, respectively. For the molecular surface they are
$\approx .9$ and $\approx .4$.

\FloatBarrier 

\section{Proof that the three extrapolation equations can be solved }

We prove that the system of three equations \eqref{front} or \eqref{form5} can always be solved provided
\beq  0 < \rho_1 < \rho_2 < \rho_3   \eeq
To do this we show that the $3 \times 3$ determinant $D$ whose $i$th row is
\beq  [ 1 \quad \rho_i I_0(x/\rho_i) \quad \rho_i^3 I_2(x/\rho_i) ] \eeq
is positive, where $x = b/h$.  For $x = 0$, the case of evaluation on the surface,
we see directly that
\beq D = (\rho_3 - \rho_2)(\rho_2 - \rho_1)(\rho_3 - \rho_1)(\rho_1 + \rho_2 + \rho_3)/(3\pi) \eeq
In general we can assume $x \geq 0$ since $I_0$ and $I_2$
are even in $x$.  First we note from \eqref{eye0} and \eqref{eye2} that
\beq I_2(x) = -\frac23 x^2 I_0(x) + \frac{1}{3\sqrt{\pi}} e^{-x^2}  \eeq
Inserting this expression in last entry of the $i$th row we obtain
\beq - \frac23 \rho_i x^2 I_0(x/\rho_i) + \frac{1}{3\sqrt{\pi}} \rho_i^3 e^{-x^2/\rho_i^2} \eeq
The third column is now a sum where the first term is a multiple of the second column.
This first part contributes zero, and the determinant becomes
\beq D = \frac{1}{3\sqrt{\pi}} \left| \begin{array}{ccc}
    1 & \rho_1 I_0(x/\rho_1) & \rho_1^3 e^{-x^2/\rho_1^2} \\
    1 & \rho_2 I_0(x/\rho_2) & \rho_2^3 e^{-x^2/\rho_2^2} \\
    1 & \rho_3 I_0(x/\rho_3) & \rho_3^3 e^{-x^2/\rho_3^2} 
\end{array} \right|   \eeq

Next we subtract row 1 from rows 2 and 3, resulting in the $2\times 2$ determinant
\beq 3\sqrt{\pi} D = \left| \begin{array}{cc}
   \rho_2 I_0(x/\rho_2) - \rho_1 I_0(x/\rho_1) \;\;\; & \rho_2^3 e^{-x^2/\rho_2^2} - \rho_1^3 e^{-x^2/\rho_1^2} \\
   \rho_3 I_0(x/\rho_3) - \rho_1 I_0(x/\rho_1) \;\;\; & \rho_3^3 e^{-x^2/\rho_3^2} - \rho_1^3 e^{-x^2/\rho_1^2}
\end{array} \right|   \eeq
We can assume that $\rho_1 = 1$,
since we could replace arbitrary $x$ with $x/\rho_1$. The new determinant has the form
\beq \label{detratio}
A(x,\rho_2)B(x,\rho_3) - A(x,\rho_3)B(x,\rho_2) = B(x,\rho_2)B(x,\rho_3)
  \left(\frac{A(x,\rho_2)}{B(x,\rho_2)} - \frac{A(x,\rho_3)}{B(x,\rho_3)} \right) \eeq
where
\beq A(x,r) = rI_0(x/r) - I_0(x) \,,\quad  B(x,r) = r^3e^{-x^2/r^2} - e^{-x^2}  \eeq
Clearly $B(x,r) > 0$ for $r>1$.  For $r = 1$, $A \equiv 0$ and $B \equiv 0$ and as
$r \to 1$ from above, $A(x,r)/B(x,r) \to A'(x,1)/B'(x,1) = 1/\left(\sqrt{\pi}(2x^2 + 3)\right)$,
as seen from \eqref{derivs} below.  Hereafter $'$ means $\pa/\pa r$.
To show $D > 0$ it suffices, according to \eqref{detratio}, to show that
$A(x,r)/B(x,r)$ decreases as $r > 1$ increases.
To verify this we will show that $(A/B)' < 0$ or equivalently
\beq \label{Fclaim}
  F(x,r) \equiv e^{x^2/r^2}\left( A(x,r)B'(x,r) - A'(x,r)B(x,r) \right) > 0 
           \,,\quad r > 1 \eeq
At $r=1$, $F \equiv 0$ since $A = B = 0$.  We find after some cancellation that
\beq \label{derivs}
e^{x^2/r^2} A' = 1/\sqrt{\pi}\,, \quad e^{x^2/r^2} B' = 2x^2 + 3r^2\,, \quad F' = 6rA\eeq
Then $A>0$ for $r>1$, since $A' > 0$ and $A(x,1) = 0$.
Finally $F' = 6rA > 0$, and since $F(x,1) = 0$, we conclude that 
$F(x,r) > 0$ for $r > 1$, as claimed in \eqref{Fclaim}.


\section{Conclusions and future work}

We have developed a simple, self-contained method for computing surface integrals,
such as single or double layer integrals for harmonic functions or for Stokes flow,
when evaluated at points close to the surface, so that these integrals are nearly singular. The integral kernel is replaced by a regularized form.  The modified integral expressions are given in Sect.~2. Asymptotic analysis in Sect.~3 provides a formula for the leading error due to this regularization, uniform for target points near the surface.  This formula can
be used with extrapolation to obtain high order regularization.  The high order allows
the modified integrands to be smooth enough so that a conventional quadrature can be used;
see Sect.~4.  Numerical tests in Sect.~5 verify the accuracy by evaluating known solutions
at points near the surface.

The tests in this work used direct summation so that errors are measured unambiguously. To reduce the high computational cost for large systems, fast summation methods such as treecodes or fast multipole methods can be used. In the present work, the integrals are computed for several values of the regularization parameter $\delta$ to obtain the extrapolated value. Since the contribution from $\delta$ decays rapidly away from the near singularity, the evaluation of the integrals for additional $\delta$ values might need to be done only in a certain neighborhood of the target point. This approach will be investigated in future work.

Surface integrals considered here are nearly singular when values are needed at grid points near the surface, which was the focus of our tests. The near singularity also occurs when multiple surfaces are close to each other. One such example was presented in our earlier work with Stokes surface integrals~\cite{tbjcp}, where corrections were added to improve the accuracy. The correction formulas are found using asymptotic analysis somewhat similar to the analysis presented here, and they have complicated expressions. Furthermore, the corrections improve the accuracy only to $O(h^3)$ in the nearly singular case. The extrapolation method presented here is more accurate and much easier to use. We therefore expect the current method to work better in multi-surface cases.

This method could be used to simulate moving interfaces in Stokes flow.  A possible approach is to represent the surface by a level set function.  To move the surface, the current velocity can be computed at grid points nearby, then the level set function is updated at these grid points, and finally the new surface is recovered.  The method developed in this work is well suited to find the velocity at the grid points, and its simplicity should be an advantage.

\section*{Declarations}

\subsection*{Conflict of interest}

The authors declare no competing interests.

\section*{Acknowledgment}
The work of ST was supported by the National Science Foundation grant
DMS-2012371.

\bibliographystyle{plain}
\bibliography{Bibtrap2}

\end{document}